%% file: main.tex
\documentclass[letterpaper,11pt,twoside]{article}
\usepackage{amsmath,amssymb,amsthm}
\usepackage{epic,eepic}
\usepackage{fancyhdr}
\usepackage{verbatim}
\usepackage{appendix}
\usepackage{latexsym}
\usepackage{enumerate}
\usepackage{mathrsfs}

\input{macrodef}

\allowdisplaybreaks

\title{\textbf{Generating connected and biconnected graphs}}
\author{\^Angela Mestre\footnote{email: mestrang@kmlinux.fjfi.cvut.cz}
\\\\
Doppler Institute \&
Department of Mathematics,\\
FNSPE,  Czech Technical University\\
Trojanova 13, 120 00 Praha 2\\
Czech Republic\\\\
Centro de F\'{\i}sica Te\'orica,\\ Departamento de
F\'{\i}sica da Universidade de Coimbra\\ P-3004-516 Coimbra,
Portugal\\
}
\begin{document}

\maketitle

\begin{abstract}
\input{abstract}

\end{abstract}

\input{intro}

\input{basics}

\input{linear}

\input{loops}

\input{biconn}

\input{simple}

\input{noloops}

\input{theta}

\input{algo}

\input{acknowledge}
\appendix
\input{apploops}
\input{appbiconn}
\input{appsimple}
\input{appnoloops}
\bibliographystyle{amsplain}
\bibliography{bibliography}
\end{document}

%% file: macrodef.tex
\newcommand{\two}{2}
\newcommand{\defeq}{:=}
\newtheorem{prop}{Proposition}
\newtheorem{thm}[prop]{Theorem}

\newtheorem{lem}[prop]{Lemma}

\newcommand{\Q}{\mathbb{Q}}
\newcommand{\card}{\mbox{card}}
\newcommand{\graph}{G}
\newcommand{\legs}{e}
\newcommand{\Vs}{V_{{\tiny\mbox{simple}}}}
\newcommand{\Vb}{V_{{\tiny\mbox{biconn}}}}
\newcommand{\Vc}{V_{{\tiny\mbox{conn}}}}
\newcommand{\Vm}{V_{{\tiny\mbox{loopless}}}}
\newcommand{\Vd}{V_{{\tiny\mbox{disconn}}}}
\newcommand{\multiedges}{\rho}
\newcommand{\lp}{k}
\newcommand{\vertex}{n}
\newcommand{\edge}{m}
\newcommand{\ext}{s}
\newcommand{\Ti}{t_i}
\newcommand{\Qci}{q^{c{(\rho)}}_{i}}
\newcommand{\Qdi}{q^{d{(\rho)}}_{i}}
\newcommand{\Qcione}{q^{c(1)}_{i}}
\newcommand{\Qdione}{q^{d(1)}_{i}}
\newcommand{\Qi}{q_{i}^{(\rho)}}
\newcommand{\Qione}{q_{i}^{(1)}}
\newcommand{\Ebij}{{l^b}_{i,j}}
\newcommand{\Eaij}{{l^a}_{i,j}}
\newcommand{\Eij}{l_{i,j}}
\newcommand{\Ein}{l_{i,n+1}}
\newcommand{\Einn}{l_{i,n}}
\newcommand{\si}{s_{\mathcal{E}_i}}
\newcommand{\sci}{{s^c}_{\mathcal{E}_i}}
\newcommand{\sdi}{{s^d}_{\mathcal{E}_i}}
\newcommand{\cij}{c_{i,j}}
\newcommand{\cin}{c_{i,n}}
\newcommand{\Ext}{E_{{\tiny\mbox{ext}}}}
\newcommand{\Int}{E_{{\tiny\mbox{int}}}}
\newcommand{\mapi}{\varphi_{{\tiny\mbox{int}}}}
\newcommand{\mape}{\varphi_{{\tiny\mbox{ext}}}}

\newcommand{\legsii}{\mathcal{E}_{{\tiny\mbox{int}},i}}
\newcommand{\legsij}{\mathcal{E}_{{\tiny\mbox{int}},j}}
\newcommand{\legsie}{\mathcal{E}_{{\tiny\mbox{ext}},i}}

\newcommand{\edgesii}{\mathcal{L}_{{\tiny\mbox{int}},i}}
\newcommand{\edgesie}{\mathcal{L}_{{\tiny\mbox{ext}},i}}
\newcommand{\edgesje}{\mathcal{L}_{{\tiny\mbox{ext}},j}}
\newcommand{\edgesij}{\mathcal{L}_{{\tiny\mbox{int}},j}}

\newcommand{\mapeii}{\varphi_{{\tiny\mbox{ext}}}^*}
\newcommand{\aux}{\eta}
\newcommand{\auxz}{\zeta}
\newcommand{\auxii}{\eta^*}
\newcommand{\auxzii}{\zeta^*}
\newcommand{\Dl}{\omega}
\newcommand{\Ds}{\theta}
\newcommand{\Dp}{\beta}
\newcommand{\Dnew}{\sigma}
\newcommand{\Dsnew}{\hat{\theta}}
\newcommand{\exta}{\epsilon}
\newcommand{\extb}{\xi}
\newcommand{\A}{\mathscr{A}}
\newcommand{\B}{\mathscr{B}}
\newcommand{\C}{\mathscr{C}}
\newcommand{\D}{\mathscr{D}}
\newcommand{\CC}{\mathcal{C}}

%% file: abstract.tex
We focus on the algorithm underlying the main result of \cite{MeOe:loop}. 
This is an algebraic formula to generate all connected graphs in a recursive and efficient manner. The key feature 
is that each graph carries a scalar factor given by the inverse of the order of its group of automorphisms.
 In the present paper, we  revise that algorithm on the level of  graphs. Moreover, we extend the result subsequently to further classes of connected graphs, namely, (edge) biconnected,  simple and loopless graphs. Our method consists of basic graph transformations only.

%% file: intro.tex
 The present paper is part of a program laid out in \cite{MeOe:npoint, MeOe:loop} with the focus on the combinatorics of different kinds of connected graphs and  problems of graph generation. In particular,
the
 main result of \cite{MeOe:npoint} is a recursion formula to generate all  trees.   This  result  is generalized  to all
 connected graphs  in  \cite{MeOe:loop}. The underlying structure is  a Hopf algebraic representation of graphs. In both cases, in a recursion step, the formulas yield  linear combinations of graphs with rational coefficients.  The essential property is that
 the coefficients of graphs are given by the inverses of the orders of
their groups of automorphisms. Other problems in this context are  considered in \cite{Jordan,Polya}, for instance.

In this paper, we express the algebraic recursion  formula to generate all connected graphs given in \cite{MeOe:loop}, in terms of  graphs. Moreover,  we extend this result successively to (edge) biconnected,  simple and loopless (connected) graphs.  Crucially, as in \cite{MeOe:npoint, MeOe:loop}, the exact coefficients of graphs are obtained. 

Our method is based on three linear graph transformations to  produce a graph  with, say, $\edge$  edges from a graph with $\edge-1$ edges. Namely,  (a) assigning a loop to a vertex; (b) connecting a pair of vertices with an edge; (c) splitting a vertex in two, distributing the ends of edges assigned to the split vertex, between the two new ones in a given way, and connecting the two new vertices with an  edge. In particular, the last operation is  (equivalently) defined   for simple graphs in \cite{103}.

Furthermore, we consider a  definition of graph which is more general than the one given in most textbooks on graph theory. In particular,  we allow edges not to be connected to vertices at both ends. Clearly, all results hold when the number of these \emph{external} edges vanishes and the standard definition of graph is recovered. However, as in \cite{MeOe:npoint, MeOe:loop}, external edges are fundamental for the (induction) proofs. This is due to the fact that vertices carrying (labeled) external edges are distinguishable and thus held fixed under any symmetry.

This paper is organized as follows: Section~\ref{sec:basics} reviews the basic concepts of graph theory that underly much of the paper. Section~\ref{sec:linear} contains the definitions of the basic linear maps to be  used in  the following sections. Section \ref{sec:loops}  translates the recursion formula to generate all connected graphs given in \cite{MeOe:loop}, to the language of graph theory. Section \ref{sec:extensions} extends this result to  biconnected, simple and loopless (connected) graphs. 
The appendixes list all
connected graphs (up to three edges), all biconnected graphs (up to four edges),  all simple connected graphs (up to five  edges), and all loopless connected graphs (up to four  edges),
with no external edges and together with their scalar factors.

%% file: basics.tex
\section{Graphs}
\label{sec:basics}
We  briefly review  the basic concepts of graph theory that are relevant for the following sections.
For more information on these we refer the reader to standard textbooks such as \cite{diestel}.

\bigskip

Let  $A$ and $B$ denote sets.  By  $[A,B]$, we denote the set of all unordered pairs of elements of $A$ and $B$, $\{a\in A, b\in B\}$. In particular, by  $[A]^\two\defeq[A,A]$, we denote the set of all $\two$-element subsets of $A$.   Also, by $\two^A$, we denote the power set of $A$, i.e., the set of all subsets of $A$.  By $\card(A)$, we denote the cardinality of the set $A$. Finally, we recall  that  the symmetric difference of the sets $A$ and $B$ is given by $A\triangle B\defeq(A\cup B)\backslash (A\cap B)$.

\medskip

Let  $\mathcal{V}=\{v_i\}_{i\in\mathbb{N}}$ and $\mathcal{K}=\{\legs_a\}_{a\in\mathbb{N}}$ be infinite sets so that $\mathcal{V}\cap\mathcal{K}=\emptyset$. Let $V\subset\mathcal{V}$; $V\neq\emptyset$ and $K\subset\mathcal{K}$ be finite
sets. 
Let  $E=\Int\cup\Ext\subseteq[K]^\two$ and $\Int\cap\Ext=\emptyset$. Also, let the elements of $E$ satisfy
$\{\legs_a,\legs_{a'}\}\cap\{\legs_b,\legs_{b'}\}=\emptyset$. 
In this context, a  \emph{graph} is a triple $\graph=(V,K,E)$ together with  the following  maps: 
\begin{enumerate}[(a)]
\item$\mapi\defeq\aux\circ\auxz:\Int\rightarrow [V]^\two\cup V; \{\legs_a,\legs_{a'}\}\mapsto \{v_i,v_{i'}\}$, 
where
\begin{itemize}
\item $\auxz:\Int\to[K,V]\cup[K]^\two$; $\auxz(\{\legs_a,\legs_{a'}\})=\{\legs_a,v_{i'}\}\quad \mbox{or}\quad\auxz(\{\legs_a,\legs_{a'}\})=\{\legs_a,\legs_{a'}\}$; 
\item $\aux:[K,V]\cup[K]^\two\to[V]^\two\cup V$; $\aux(\{\legs_a,v_{i'}\})=\aux(\{\legs_a,\legs_{a'}\})=\{v_{i},v_{i'}\}$;  \end{itemize} 
\item $\mape:\Ext\rightarrow [V,K]; \{\legs_a,\legs_{a'}\}\mapsto \{v_i,\legs_{a'}\}$. 
\end{enumerate}
The elements of $V$ and $E$ are called \emph{vertices} and \emph{edges}, respectively. In particular, the elements of $\Int$ and $\Ext$ are called  \emph{internal} edges and \emph{external} edges, respectively. Both internal and external edges correspond to unordered pairs of elements of $K$. The elements of  these pairs are called \emph{ends} of edges.
In other words, internal edges are edges that are connected to vertices at both ends, while external  edges have one free end. Internal edges with both ends assigned to the
same vertex are also called \emph{loops}.
Two distinct vertices  connected together by one or more internal edges, are said to be \emph{adjacent}. Two or more internal edges connecting the same pair of distinct vertices together, are called \emph{multiple edges}. For instance, Figure  \ref{fig:3graphs} (a)  shows a loop, while Figure \ref{fig:3graphs} (b) shows  a graph with both multiple edges and  external edges.
\begin{figure}
\begin{center}
\input{figures/3graphs.eepic}\label{fig:3graphs}
\caption{ (a) A loop; (b) A graph with both internal and external edges; (c) A graph with labeled external edges.  Internal edges are  represented by continuous lines, while the external ones are represented by  dashed lines.}
\end{center}
\end{figure}
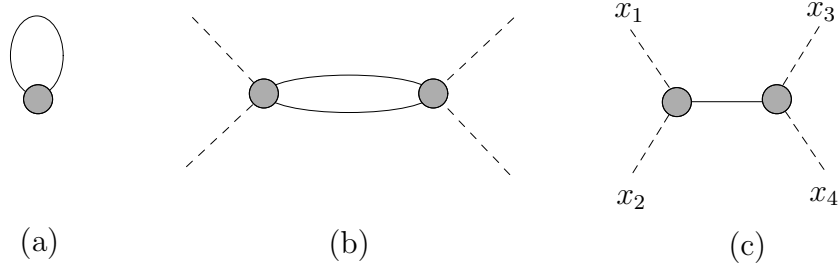  A graph with no loops nor multiple edges is called \emph{simple}. The \emph{degree} of a vertex is the number of ends of edges
assigned to the vertex.

Let $\graph=(V,K,E)$; $E=\Int\cup\Ext$,  $\card(\Ext)=\ext$, together with the maps $\mapi$ and $\mape$, denote a graph.  
The external edges of the graph $\graph$ are said to be \emph{labeled} if their  free ends are assigned labels
$x_1,\dots, x_{\ext}$ from a  \emph{label set} $L=\{x_1,\ldots,x_{\ext}\}$. Labels on different ends of
external edges are required to be distinct. In other words, a \emph{labeling} of the external edges of the graph $\graph$, is an injective map $l:\Ext\rightarrow [K,L]; \{\legs_a,\legs_{a'}\}\mapsto \{\legs_a,x_z\}$, where $z\in\{1,\ldots,\ext\}$.  For instance, Figure  \ref{fig:3graphs} (c)  shows a graph  with two vertices and four labeled external edges. 

\medskip

A graph
$\graph^*=(V^*,K^*,E^*)$; $E^*=\Int^*\cup\Ext^*$, together with the maps $\mapi^*$ and $\mape^*$, is called a \emph{subgraph} of a graph
$\graph=(V,K,E)$; $E=\Int\cup\Ext$, together with the maps $\mapi$ and $\mape$, if $V^*\subseteq V$, $K^*\subseteq K$,
$E^*\subseteq E$ and $\mapi^*=\mapi\rvert_{\Int^*}$, $\mape^*=\mape\rvert_{\Ext^*}$.

A \emph{path}  is  a graph $P=(V,K,\Int)$; $V=\{v_1,\ldots,v_{\vertex}\}$, $\vertex\defeq\card(V)>1$, together with the map $\mapi$,
so that   $\mapi(\Int)=\{\{v_1,v_\two\},\{v_\two, v_3\},\ldots,\{v_{n-1},v_\vertex\}\}$ and the vertices $v_1$ and $v_\vertex$ have degree 1, while the vertices $v_\two,\ldots,v_{\vertex-1}$ have degree $\two$. In this context, the vertices $v_1$ and $v_\vertex$ are called the \emph{end
 point} vertices, while the vertices $v_\two,\ldots,v_{\vertex-1}$ are called the \emph{inner} vertices.
A \emph{cycle}   is a  graph $C=(V',K',\Int')$; $V'=\{v_1,\ldots,v_{\vertex}\}$,  together with the map $\mapi'$,
 so that   $\mapi'(\Int')=\{\{v_1,v_\two\},\{v_\two, v_3\},\ldots,\{v_{\vertex-1},v_\vertex\}, \{v_\vertex, v_1\}\}$ and every vertex has degree $\two$. 
A graph is said to be \emph{connected} if every pair of vertices is joined by a path. Otherwise, it is \emph{disconnected}.
Moreover, a {\em tree}  is a connected graph with no cycles.
 A {\em biconnected} graph (or \emph{edge-biconnected} graph) is a connected graph that remains connected after erasing one and whichever internal edge. 
By definition,  a graph consisting of a single vertex  is  biconnected. 

Furthermore, let  $\graph=(V,K,E)$; $E=\Int\cup\Ext$, together with the maps $\mapi$ and $\mape$ denote a graph. The set $\two^{\Int}$ is a vector space over the field $\mathbb{Z}_\two$ so that vector addition is given by the  symmetric difference. The \emph{cycle space} $\CC$ of the graph $\graph$ is defined as the subspace of $\two^{\Int}$ generated by all the cycles of $\graph$. The dimension of $\CC$ is  called the \emph{cyclomatic number} of the graph $\graph$.  Moreover, the cyclomatic number $\lp\defeq\dim\CC$  yields in terms of the  vertex number $\vertex\defeq\mbox{card}(V)$ and the internal edge number $\edge\defeq\mbox{card}(\Int)$ as $\lp=\edge-\vertex+c$, where $c$ denotes the number of connected components of the graph $\graph$ \cite{Kirch}.

\medskip

Now, let $L=\{x_1,\ldots,x_{\ext}\}$ be a finite label set. Let $\graph=(V,K,E)$; $E=\Int\cup\Ext$, together with the maps $\mapi$ and $\mape$, and $\graph^*=(V^*,K^*,E^*)$;  $E^*=\Int^*\cup\Ext^*$, together with the maps $\mapi^*$ and $\mape^*$, denote two graphs. Let $l:\Ext\rightarrow [K,L]$ and $l^*:\Ext^*\rightarrow [K^*,L]$ be labelings of the elements of $\Ext$ and of $\Ext^*$, respectively.  An \emph{isomorphism} between the graphs $\graph$ and $\graph^*$ is a bijection $\psi_V:V\to V^*$ and  a bijection $\psi_K:K\to K^*$  which satisfy the following three conditions:
\begin{enumerate}[(a)]
\item If  $\mapi(\{\legs_{a},\legs_{a'}\})=\{v_i,v_{i'}\}$ then $\mapi^*(\{\psi_K(\legs_{a}),\psi_K(\legs_{a'})\})=\{\psi_V(v_i),\psi_V(v_{i'})\}$; 
\item If  $\mape(\{\legs_{a},\legs_{a'}\})=\{v_i,\legs_{a'}\}$ then $\mape^*\{\psi_K(\legs_{a}),\psi_K(\legs_{a'})\}=\{\psi_V(v_{i}),\psi_K(\legs_{a'}) \}$; 
\item $L\cap l(\{\legs_a,\legs_{a'}\})=L\cap l^*(\{\psi_K(\legs_a),\psi_K(\legs_{a'})\})$.
\end{enumerate}
Clearly, an isomorphism  defines an equivalence relation on graphs.
In particular, a \emph{vertex} (resp. \emph{edge}) isomorphism between the graphs $\graph$ and $\graph^*$  is an isomorphism
so that $\psi_E$ (resp. $\psi_V$) is the identity map. 
In this context, a \emph{symmetry} of a graph $\graph$, is an isomorphism
 of the graph  onto itself (i.e, an \emph{automorphism}). The order of the group  of automorphisms of the graph $\graph$  is called the \emph{symmetry
factor},  denoted by $S^\graph$.
 Also,
a \emph{vertex symmetry} (resp. \emph{edge symmetry}) of
a   graph $\graph$, is a vertex (resp. edge) automorphism of the graph.  The
order of the group of vertex (resp. edge) automorphisms is called the \emph{vertex
  symmetry factor} (resp. \emph{edge
  symmetry factor}) of the graph.  This is denoted  by
$S^\graph_{\text{vertex}}$ (resp. $S^\graph_{\text{edge}}$).
Furthermore,  the orders of the groups of vertex and edge automorphisms of a graph $\graph$,  satisfy
$S^{\graph}=
S^{\graph}_{\text{vertex}}\cdot
S^{\graph}_{\text{edge}}$ (a proof is given in \cite{MeOe:loop}, for instance).

%% file: figures/3graphs.eepic
\setlength{\unitlength}{0.00083333in}
\begingroup\makeatletter\ifx\SetFigFontNFSS\undefined%
\gdef\SetFigFontNFSS#1#2#3#4#5{%
  \reset@font\fontsize{#1}{#2pt}%
  \fontfamily{#3}\fontseries{#4}\fontshape{#5}%
  \selectfont}%
\fi\endgroup%
{\renewcommand{\dashlinestretch}{30}
\begin{picture}(5096,1758)(0,-10)
\put(173,1315){\ellipse{330}{496}}
\put(2121,1075){\ellipse{1042}{234}}
\dashline{60.000}(1108,621)(1596,1078)
\dashline{60.000}(2661,1092)(3149,1549)
\dashline{60.000}(2668,1057)(3125,569)
\dashline{60.000}(1146,1553)(1603,1065)
\put(67,93){\makebox(0,0)[lb]{\smash{{\SetFigFontNFSS{12}{14.4}{\familydefault}{\mddefault}{\updefault}(a)}}}}
\put(1994,83){\makebox(0,0)[lb]{\smash{{\SetFigFontNFSS{12}{14.4}{\familydefault}{\mddefault}{\updefault}(b)}}}}
\dashline{60.000}(3878,1475)(4172,1035)
\dashline{60.000}(4790,1043)(5084,603)
\dashline{60.000}(4780,1055)(5053,1495)
\dashline{60.000}(3893,583)(4166,1023)
\path(4205,1027)(4779,1027)
\texture{555555 55000000 555555 55888888 88555555 55000000 555555 55808080 
	80555555 55000000 555555 55888888 88555555 55000000 555555 55888088 
	80555555 55000000 555555 55888888 88555555 55000000 555555 55808080 
	80555555 55000000 555555 55888888 88555555 55000000 555555 55888088 }
\put(182,1041){\shade\ellipse{182}{182}}
\put(182,1041){\ellipse{182}{182}}
\put(4166,1025){\shade\ellipse{182}{182}}
\put(4166,1025){\ellipse{182}{182}}
\put(4790,1042){\shade\ellipse{182}{182}}
\put(4790,1042){\ellipse{182}{182}}
\put(1590,1078){\shade\ellipse{182}{182}}
\put(1590,1078){\ellipse{182}{182}}
\put(2647,1078){\shade\ellipse{182}{182}}
\put(2647,1078){\ellipse{182}{182}}
\put(4488,73){\makebox(0,0)[lb]{\smash{{\SetFigFontNFSS{12}{14.4}{\familydefault}{\mddefault}{\updefault}(c)}}}}
\put(3777,1553){\makebox(0,0)[lb]{\smash{{\SetFigFontNFSS{12}{14.4}{\rmdefault}{\mddefault}{\updefault}$x_1$}}}}
\put(4974,1560){\makebox(0,0)[lb]{\smash{{\SetFigFontNFSS{12}{14.4}{\rmdefault}{\mddefault}{\updefault}$x_3$}}}}
\put(4993,408){\makebox(0,0)[lb]{\smash{{\SetFigFontNFSS{12}{14.4}{\rmdefault}{\mddefault}{\updefault}$x_4$}}}}
\put(3789,386){\makebox(0,0)[lb]{\smash{{\SetFigFontNFSS{12}{14.4}{\rmdefault}{\mddefault}{\updefault}$x_2$}}}}
\end{picture}
}

%% file: linear.tex
\section{Elementary linear transformations}\label{sec:linear}
We introduce some linear maps and prove their fundamental properties.

\bigskip

Given an arbitrary set $W$, by $\Q W$, we  denote the free vector space on the  set  $W$ over $\Q$.  That is,
 (a) every vector in $\Q W$ yields a linear combination of the elements of $W$ with coefficients in $\Q$;  (b) the set $W$ is linearly independent. 

\medskip

Let  $\mathcal{V}=\{v_i\}_{i\in\mathbb{N}}$ and $\mathcal{K}=\{\legs_a\}_{a\in\mathbb{N}}$ be infinite sets so that $\mathcal{V}\cap\mathcal{K}=\emptyset$. Fix integers $t,\ext,\lp\ge0$  and  $\vertex\ge1$. Let $L=\{x_1,\ldots,x_\ext\}$ be a   label set. By $V^{\vertex,\lp,\ext}$, we denote the set  of all
graphs with $\vertex$ vertices,  cyclomatic number $\lp$ and $\ext$ external edges whose free ends are labeled $x_1,\ldots,x_\ext$.  In all that follows, let $V=\{v_1,\ldots,v_\vertex\}\subset\mathcal{V}$, $K=\{\legs_1,\ldots,\legs_t\}\subset\mathcal{K}$ and $E=\Int\cup\Ext$ be the sets of vertices, of ends of edges and of edges, respectively, of all elements of $V^{\vertex,\lp,\ext}$. 
Also,  let $l:\Ext\rightarrow [K,L]$ be a labeling of  their external edges.  Moreover, by $\Vc^{\vertex,\lp,\ext}$ and $\Vd^{\vertex,\lp,\ext}$,   we denote the subsets of $V^{\vertex,\lp,\ext}$ whose elements are connected or disconnected graphs, respectively. Finally, by $\Vb^{\vertex,\lp,\ext}$, $\Vs^{\vertex,\lp,\ext}$ and $\Vm^{\vertex,\lp,\ext}$,  we denote the subsets of $\Vc^{\vertex,\lp,\ext}$ whose elements are biconnected,  simple and loopless graphs, respectively. 

We now define the following  linear transformations:

\begin{enumerate}[(i)]
\item
\emph{Assigning a loop to a vertex}: Let $\graph=(V,K,E)$ 
together with the maps $\mapi$ and $\mape$, denote a graph in $V^{\vertex, \lp,\ext}$. For all $i\in\{1,\ldots,\vertex\}$, define 
$$\Ti:\Q V^{\vertex,\lp,\ext}\rightarrow\Q V^{\vertex,\lp+1,\ext};\graph\mapsto\graph^*\,,$$ 
where  the graph $\graph^*=(V^*,K^*,E^*)$; $E^*=\Int^*\cup\Ext^*$, together with the maps $\mapi^*$ and $\mape^*$, satisfies the following conditions:
\begin{enumerate}[(a)]
\item $V^*=V$;\item $K^*=K\cup\{\legs_{t+1},\legs_{t+\two}\}$; 
\item $E^*=\Int^*\cup\Ext^*$, where $\Int^*=\Int\cup\{\legs_{t+1},\legs_{t+\two}\}$, $\Ext^*=\Ext$;  
\item $\mapi^*\rvert_{\Int}=\mapi$ and $\mapi^*(\{\legs_{t+1},\legs_{t+\two}\})=\{v_i\}$; 
\item $\mape^*=\mape$.
\end{enumerate} 
The $\Ti$-maps
are extended to all of $\mathbb{Q} V^{\vertex,\lp,\ext}$ by linearity.
Since the map $\Ti:\Q V^{\vertex,\lp,\ext}\rightarrow\Ti(\Q V^{\vertex,\lp,\ext})$ is injective,  the operation of \emph{erasing a loop}  is given by  $\Ti^{-1}$. 

\item
\emph{Connecting a pair of distinct vertices with an internal edge}: Let $\graph=(V,K,E)$ 
together with the maps $\mapi$ and $\mape$, denote a graph in $V^{\vertex, \lp,\ext}$; $\vertex>1$. For  all $i,j\in\{1,\ldots,\vertex\}$ with $i\neq j$,  define 
$$\Eij:\Q V^{\vertex,\lp,\ext}\rightarrow\Q V^{\vertex,\lp+1,\ext}\cup\Q V^{\vertex,\lp,\ext};\graph\mapsto\graph^*\,,$$
where  the graph $\graph^*=(V^*,K^*,E^*)$; $E^*=\Int^*\cup\Ext^*$, together with the maps $\mapi^*$ and $\mape^*$, satisfies the following conditions:
\begin{enumerate}[(a)]
\item $V^*=V$;
\item $K^*=K\cup\{\legs_{t+1},\legs_{t+\two}\}$; 
\item $E^*=\Int^*\cup\Ext^*$, where $\Int^*=\Int\cup\{\legs_{t+1},\legs_{t+\two}\}$, $\Ext^*=\Ext$;  
\item $\mapi^*\rvert_{\Int}=\mapi$ and $\mapi^*(\{\legs_{t+1},\legs_{t+\two}\})=\{v_i,v_j\}$; 
\item $\mape^*=\mape$.
\end{enumerate} 
The $\Eij$-maps
are extended to all of $\mathbb{Q}V^{\vertex,\lp,\ext}$ by linearity.
Since the map $\Eij:\Q V^{\vertex,\lp,\ext}\rightarrow\Eij(\Q V^{\vertex,\lp,\ext})$ is injective,   the operation of \emph{erasing an internal edge distinct from a loop} is given by   $\Eij^{-1}$.
Furthermore, for $\vertex>1$ (resp. $\vertex>\two$) and for all $i,j\in\{1,\ldots,\vertex\}$ with $i\neq j$,  define $\Eaij\defeq\Eij\circ\delta_{i,j}$ (resp. $\Ebij\defeq\Eij\circ(\mbox{id}-\delta_{i,j}$)), where $\mbox{id}:\Q V^{\vertex,\lp,\ext}\to\Q V^{\vertex,\lp,\ext}$ is the identity map and $\delta_{i,j}:\Q V^{\vertex,\lp,\ext}\to\Q V^{\vertex,\lp,\ext};\graph\mapsto\left\{\begin{array}{ccc}\graph & \mbox{if} & \{v_i,v_j\}\in\mapi(\Int)\\0 & \mbox{otherwise} & \\\end{array}\right.$ is a linear map. 

\item
\emph{Splitting a vertex in two and distributing the ends of edges assigned to the split vertex,
between the two new ones in all possible ways}:
 Let $\graph=(V,K,E)$ 
together with the maps $\mapi\defeq\aux\circ\auxz$ and $\mape$, denote a graph in $V^{\vertex,\lp,\ext}$. 
 Let $\mathcal{L}_i\subseteq E$ be the set of edges connected to the vertex $v_i\in V$; $i\in\{1,\ldots, \vertex\}$. Also, let $\edgesii$ and $\edgesie$ be the subsets of $\mathcal{L}_i$ whose elements are internal edges or external edges, respectively.  Hence, $\mathcal{L}_i=\edgesii\cup\edgesie$ and $\edgesii\cap\edgesie=\emptyset$.  Moreover, let   $\mathcal{E}_i\subseteq K$ be  the set of ends of edges assigned to the vertex $v_i$. Also, let $\legsii$ and $\legsie$ be the subsets of $\mathcal{E}_i$ whose elements are  ends of internal edges or ends of external edges, respectively. Thus, $\mathcal{E}_i=\legsii\cup\legsie$ and $\legsii\cap\legsie=\emptyset$. Let $[\legsii']^ \two\defeq\edgesii\cap[\legsii]^\two$ and $\legsii''\defeq\legsii\backslash\legsii'$. Finally,  let $\mathcal{I}^\two_{\mathcal{E}_i}$ denote the set of all partitions of the set $\mathcal{E}_i$  into two disjoint sets: $\mathcal{I}^\two_{\mathcal{E}_i}=\{\{\mathcal{E}_i^{(1)},\mathcal{E}_i^{(\two)}\}:\mathcal{E}_i^{(1)}\cup \mathcal{E}_i^{(\two)}=\mathcal{E}_i\quad  \mbox{and}\quad \mathcal{E}_i^{{(1)}}\cap \mathcal{E}_i^{(\two)}=\emptyset\}$. Clearly,  $\mathcal{E}_i^{(b)}= \legsii^{(b)}\cup\legsie^{(b)}$; $b\in\{1,\two\}$. Also, a partition of the set  $\legsie$, generates a partition of the set  $\edgesie$. Hence,  $\edgesie^{(b)}\subset [\legsie^{(b)}, K]$. In this context, for all $i\in\{1,\ldots,\vertex\}$, define 
$$\si:\Q V^{\vertex,\lp,\ext}\rightarrow\Q V^{\vertex+1,\lp-1,\ext}\cup\Q V^{\vertex+1,\lp,\ext};\graph\mapsto\sum_{\{\mathcal{E}_i^{(1)},\mathcal{E}_i^{(\two)}\}\in\mathcal{I}^\two_{\mathcal{E}_i}}\graph_{\{\mathcal{E}_i^{(1)},\mathcal{E}_i^{(\two)}\}}\,,$$
where the graphs $\graph_{\{\mathcal{E}_i^{(1)},\mathcal{E}_i^{(\two)}\}}=(V^*,K^*,E^*)$; $E^*=\Int^*\cup\Ext^*$, together with the maps $\mapi^*$ and $\mape^*$, satisfy the following conditions:
\begin{enumerate}[(a)]
\item $V^*=V\cup\{v_{\vertex+1}\}$; 
\item $K^*=K$;
\item $E^*=\Int^*\cup\Ext^*$, where $\Int^*=\Int$ and $\Ext^*=\Ext$;
\item  $\mapi^*=\aux^*\circ\auxz^*$, where 
\begin{itemize}
\item
$\auxzii|_{\Int\backslash\edgesii}=\auxz|_{\Int\backslash\edgesii}$; \\
 $\auxzii([\legsii']^\two)=[\legsii']^ \two$;\\ 
 $\auxzii(\edgesii\backslash[\legsii']^ \two)=[\legsii'',V']$, where $V'\subseteq V\backslash\{v_i\}$; \\
\item
$\auxii|_{\auxz(\Int\backslash\edgesii)}=\aux|_{\auxz(\Int\backslash\edgesii)}$;\\$\auxii([\legsii''^{(1)},V']\cup[\legsii'^{(1)}]^\two)=\{v_i\}$ and $\auxii([\legsii''^{(2)},V']\cup[\legsii'^{(2)}]^\two)=\{v_{\vertex+1}\}$, where $\legsii'^{(b)}\cup\legsii''^{(b)}=\legsii^{(b)}$; $b\in\{1,2\}$;\end{itemize}

\item $\mapeii|_{\Ext\backslash\edgesie}=\mape|_{\Ext\backslash\edgesie}$ and $\mapeii({\edgesie^{(1)}})=\{v_i\}$, $\mapeii({\edgesie^{(\two)}})=\{v_{\vertex+1}\}$.
\end{enumerate} 
The $\si$-maps
are extended to all of $\mathbb{Q} V^{\vertex,\lp,\ext}$ by linearity.
Moreover, we define the $\sci$-maps  (resp. $\sdi$-maps) by restricting the image of $\si$ to $\Q \Vc^{\vertex+1,\lp-1,\ext}$  (resp. $\Q \Vd^{\vertex+1,\lp,\ext}$). 

\bigskip

 Furthermore, let $\Eij^\multiedges\defeq\underbrace{
\Eij\circ\ldots\circ \Eij}_\text{$\multiedges$ times}$, where $\multiedges$ is an integer. We now combine the $\Ein^\multiedges$ and $\si$-maps to define the maps $$\Qi\defeq\frac{1}{\two(\multiedges-1)!}\Ein^\multiedges\circ\si:\mathbb{Q}V^{\vertex,\lp,\ext}\to\mathbb{Q}V^{\vertex+1,\lp+\multiedges-1,\ext}\,.$$
For $\multiedges=1$, the definition of $\Qi$  generalizes the basic  operation given in \cite{103}  to all partitions of the set $\mathcal{E}_i$ into two sets $\mathcal{E}_i^{(1)}$ and $\mathcal{E}_i^{(\two)}$, and to all vertices of the graph. Analogously,  the $\Qci$-maps (resp. $\Qdi$-maps) are given by the composition of  $\Ein^\multiedges$ with $\sci$  (resp. $\sdi$).

\item
In addition, we revise the operation of \emph{contracting an internal edge  connecting two distinct vertices, and fusing the two vertices into one} \cite{diestel}:   Let $\graph=(V,K,E)$ 
together with the maps $\mapi$ and $\mape$, denote a graph in $V^{\vertex,\lp,\ext}$; $\vertex>1$. 
Let $\{\legs_{a},\legs_{a'}\}\in \Int$ denote an internal edge connecting two distinct vertices, say, $v_i,v_j\in V$;  $i,j\in\{1,\ldots,\vertex\}$ with $i<j$. Let $\edgesij$, $\edgesje$ and $\legsij$ denote the sets of internal edges, of external edges and of ends of internal edges,  respectively, assigned to the vertex $v_j$. Let $[\legsij']^ \two\defeq\edgesij\cap[\legsij]^\two$. Finally, let $[V',v_j]\defeq\mapi(\edgesij\backslash[\legsij']^\two)\subseteq[V\backslash\{v_j\},v_j]$.
In this context, define 
$$\cij:\Q V^{\vertex,\lp,\ext}\rightarrow\Q V^{\vertex-1,\lp,\ext};\graph\mapsto\graph^*\,,$$
where the graph $\graph^*=(V^*,K^*,E^*)$; $E^*=\Int^*\cup\Ext^*$, together with the maps $\mapi^*$ and $\mape^*$, satisfies the following conditions:
\begin{enumerate}[(a)]
\item $V^*=\chi_v(V\backslash\{v_{j}\})$, where $\chi_v:v_l\mapsto\left\{\begin{array}{cc}v_l &\forall l\in\{1,\ldots,j-1\}\\ v_{l-1} &\forall l\in\{j+1,\ldots,\vertex\}\\
\end{array}\right.$ is a bijection;
\item $K^*=\chi_e(K\backslash\{\legs_{a},\legs_{a'}\})$, where $\chi_e:\legs_b\mapsto\left\{\begin{array}{cc} \legs_{b} &\forall b\in\{1,\ldots,\mbox{min}(a,a')-1\}\\\legs_{b-1} &\forall b\in\{\mbox{min}(a,a')+1,\ldots,\mbox{max}(a,a')-1\}\\\legs_{b-\two} &\forall b\in\{\mbox{max}(a,a')+1,\ldots,t\}\\\end{array}\right.$ is a bijection; 
\item $E^*=\Int^*\cup\Ext^*$, where $\Int^*=\Int\backslash\{\legs_a,\legs_{a'}\}$ and $\Ext^*=\Ext$;  
\item $\mapi^*|_{\Int^*\backslash\edgesij}=\mapi|_{\Int^*\backslash\edgesij}$;\\$\mapi^*([\legsij']^\two)=\{v_i\}$, $\mapi^*(\edgesij\backslash\{[\legsij']^\two\cup\{\legs_{a},\legs_{a'}\}\})=[V',v_i]$;
\item $\mape^*|_{\Ext\backslash\edgesje}=\mape|_{\Ext\backslash\edgesje}$ and  $\mape^*(\edgesje)=\{v_{i}\}$. 
\end{enumerate} 
The $\cij$-maps
are extended to all of $\mathbb{Q} V^{\vertex,\lp,\ext}$ by linearity.

\item \emph{Distributing external edges between all elements of a given subset of vertices in all possible ways}:
Let $\graph=(V,K,E)$  
together with the maps $\mapi$ and $\mape$, denote a graph in $V^{\vertex,\lp,\ext}$.  Let $V'=\{v_{z_1},\ldots,v_{z_{\vertex'}}\}\subseteq V$; with $1\le z_1<\dots< z_{\vertex'}\le \vertex$. 
Let $K'\subset\mathcal{K}$ be a finite set so that 
$K\cap K'=\emptyset$.  Also, let $\Ext'\subseteq [K']^\two$; $\card(\Ext')=\ext'$. Assume that the elements of $\Ext'$  satisfy  $\{\legs_a,\legs_{a'}\}\cap\{\legs_{b},\legs_{b'}\}=\emptyset$.  Let $L'=\{x_{\ext+1},\ldots,x_{\ext+\ext'}\}$ be a  label set so that $L\cap L'=\emptyset$. Also, let $l':\Ext'\rightarrow [K',L']$ be a labeling of the elements of $\Ext'$.  
 Finally, let  $\mathcal{I}^{\vertex'}_{\Ext'}$ denote the set of all partitions of the set $\Ext'$  into $\vertex'$ disjoint subsets:  $\mathcal{I}^{\vertex'}_{\Ext'}=\{\{\Ext'^{(1)},\ldots,\Ext'^{(\vertex')}\}:\Ext'^{(1)}\cup\ldots\cup \Ext'^{(\vertex')}=\Ext' \quad \mbox{and}\quad \Ext'^{{(i)}}\cap \Ext'^{(j)}=\emptyset\,,\forall i, j\in\{1,\ldots,\vertex'\}\quad\mbox{with}\quad i\neq j\}$. In this context, define 
$$\extb_{\Ext',V'} :\Q V^{\vertex,\lp,\ext}\rightarrow\Q V^{\vertex,\lp,\ext+\ext'};\graph\mapsto\sum_{\{\Ext'^{(1)},\ldots,\Ext'^{(\vertex')}\}\in \mathcal{I}^{\vertex'}_{\Ext'}}\graph_{\{\Ext'^{(1)},\ldots,\Ext'^{(\vertex')}\}}\,,$$ 
where the graphs $\graph_{\{\Ext'^{(1)},\ldots,\Ext'^{(\vertex)}\}}=(V^*,K^*,E^*)$; $E^*=\Int^*\cup\Ext^*$, together with the maps $\mapi^*$ and $\mape^*$, satisfy the following conditions:
\begin{enumerate}[(a)]
\item $V^*=V$;
\item $K^*=K\cup K'$;
\item $E^*=\Int^*\cup\Ext^*$, where $\Int^*=\Int$, $\Ext^*=\Ext\cup\Ext'$;  
\item $\mapi^*=\mapi$; 
\item $\mape^*|_{\Ext}=\mape$ and $\mape^*({\Ext'^{(i)}})=\{v_{z_i}\},\forall i\in\{1,\ldots,\vertex'\}$; 
\item $l^*:\Ext^*\rightarrow [K^*,L\cup L']$, with $l^*|_{\Ext}=l$ and $l^*|_{\Ext'}=l'$, is a labeling of the elements of $\Ext^*$.
\end{enumerate} 
The $\extb_{\Ext',V'}$-maps
are extended to all of $\Q V^{\vertex,\lp,\ext}$ by linearity.

\item \emph{Assigning external edges to vertices which have none}: Let $\graph=(V,K,E)$ 
together with the maps $\mapi$ and $\mape$, denote a graph in $V^{\vertex,\lp,\ext}$.   Assume that there exists a set $V'\subseteq V$; $\mbox{card}(V')=\ext'$ so that $V'\cap\mape(\Ext)=\emptyset$. Moreover, let $K'\subset\mathcal{K}$ be a finite set so that 
$K\cap K'=\emptyset$.  Also, let $\Ext'\subseteq [K']^\two$; $\card(\Ext')=\ext'$. Assume that the elements of $\Ext'$  satisfy  $\{\legs_a,\legs_{a'}\}\cap\{\legs_{b},\legs_{b'}\}=\emptyset$.  Let $L'=\{x_{\ext+1},\ldots,x_{\ext+\ext'}\}$ be a  label set so that $L\cap L'=\emptyset$. Finally, let $l':\Ext'\rightarrow [K',L']$ be a labeling of the elements of $\Ext'$.  
  In this context, define $$\exta_{\Ext'} :\Q V^{\vertex,\lp,\ext}\rightarrow\Q V^{\vertex,\lp,\ext+\ext'};\graph\mapsto\graph^*\,,$$ 
where  the graph $\graph^*=(V^*,K^*,E^*)$; $E^*=\Int^*\cup\Ext^*$, together with the maps $\mapi^*$ and $\mape^*$, satisfies the following conditions: 
\begin{enumerate}[(a)]
\item $V^*=V$;
\item $K^*=K\cup K'$;
\item $E^*=\Int^*\cup\Ext^*$, where $\Int^*=\Int$, $\Ext^*=\Ext\cup\Ext'$;  
\item $\mapi^*=\mapi$; 
\item $\mape^*|_{\Ext}=\mape$ and $\mape^*({\Ext'})=V'$ is a bijection;
\item $l^*:\Ext^*\rightarrow [K^*,L\cup L']$, with $l^*|_{\Ext}=l$ and $l^*|_{\Ext'}=l'$, is a labeling of the elements of $\Ext^*$.
\end{enumerate} 
The  $\exta_{\Ext'}$-maps
are extended to all of $\Q V^{\vertex,\lp,\ext}$ by linearity.
\end{enumerate}

\medskip

The following lemmas are now established.

\begin{lem}\label{lem:conntoconn} Fix integers $\lp,\ext\ge0$ and $\vertex\ge 1$. Then, for all $i,j\in\{1,\ldots,\vertex\}$ with $i\neq j$, the following statements hold:
\begin{enumerate}[(a)]\item
$\Ti(\Q \Vc^{\vertex,\lp,\ext})\subseteq\Q \Vc^{\vertex,\lp+1,\ext}$;\item
$\Eij(\Q \Vc^{\vertex,\lp,\ext})\subseteq\Q\Vc^{\vertex,\lp+1,\ext}$;\item  $\Qi(\Q\Vc^{\vertex,\lp,\ext})\subseteq\Q\Vc^{\vertex+1,\lp+\multiedges-1,\ext}$. \end{enumerate}
\end{lem}
\begin{proof}
(a), (b)  Clearly, the statements hold.
(c)  Let $\graph=(V,K,E)$ 
together with the maps $\mapi$ and $\mape$, denote a graph in $V^{\vertex, \lp,\ext}$. Let $\mathcal{E}_i$ denote the set of ends of edges assigned to the vertex $v_i\in V$.
Apply the map $\si$ to the graph $\graph$.   
In particular, $\si(\graph)$ is a linear combination of graphs, each of which is
either connected or disconnected with two components.
 Applying the $\Ein^{\multiedges}$-map to $\si(\graph)$   yields  connected graphs. This completes the proof.
\end{proof}

\begin{lem}\label{lem:onlysimple} Fix integers $\lp,\ext\ge0$ and $\vertex\ge 1$.
    Let $\graph=(V,K,E)$ 
together with the maps $\mapi$ and $\mape$, denote a graph in $V^{\vertex, \lp,\ext}$ which is not simple. Then, for all  $i, j\in\{1,\ldots,\vertex\}$ with $i\neq j$, the following statements hold:   \begin{enumerate}[(a)]\item $\Ebij(\graph)\notin\Q\Vs^{\vertex,\lp+1,\ext}$; \item $\Qdione(\graph)\notin\Q\Vs^{\vertex+1,\lp,\ext}$. \end{enumerate}
\end{lem}
\begin{proof}
 (a) Clearly, the statement holds. (b) For $\lp=0$,  the statement holds as the  $\Qdione$-maps   produce  trees   only from trees. Now, let $\vertex>1$  and $k>0$.  By assumption, the graph $\graph$  has at least either one loop or multiple edges. Let $\mathcal{E}_i$ denote the set of ends of edges assigned to the vertex $v_i\in V$. Apply the $\sdi$-map to graph $\graph$. In particular, $\sdi(\graph)$ is a linear combination of disconnected graphs, each of which is 
 produced from  the graph $\graph$ by
 assigning all ends of internal edges belonging to the same cycles, from the  vertex $v_i$ to either $v_i$ or $v_{\vertex+1}$.    Therefore, 
at least one of the two components of the graphs in $\sdi(\graph)$, is not simple. Applying the $\Ein$-map to $\sdi(\graph)$  cannot produce  simple graphs.  This completes the proof.
\end{proof}
\begin{lem}\label{lem:simpletosimple} Fix integers $\lp,\ext\ge0$ and $\vertex\ge 1$. Then, for all  $i, j\in\{1,\ldots,\vertex\}$ with $i\neq j$, the following statements hold:
\begin{enumerate}[(a)]\item $\Ebij(\Q \Vs^{\vertex,\lp,\ext})\subseteq\Q\Vs^{\vertex,\lp+1,\ext}$; \item $\Qione(\Q\Vs^{\vertex,\lp,\ext})\subseteq\Q\Vs^{\vertex+1,\lp,\ext}$.\end{enumerate}\end{lem}
\begin{proof}
 Applying the $\Ebij$ or $\Qione$-maps to a simple graph cannot introduce loops nor multiple edges.\end{proof}

\begin{lem}\label{lem:biconntobiconn}  Fix integers $\lp,\ext\ge0$ and $\vertex\ge 1$. Let $\graph=(V,K,E)$ 
together with the maps $\mapi$ and $\mape$, denote a graph in $V^{\vertex, \lp,\ext}$ which is not biconnected. Then, for all  $i\in\{1,\ldots,\vertex\}$, the following statements hold:   \begin{enumerate}[(a)]\item  $\Ti(\graph)\notin\Q\Vb^{\vertex,\lp+1,\ext}$; \item  $\Qione(\graph)\notin\Q\Vb^{\vertex+1,\lp,\ext}$.  \end{enumerate}\end{lem}
\begin{proof}
  (a) Clearly, the statement holds.
 (b) First,  by definition all connected graphs with only one vertex are biconnected. Consequently, the $\Qione$-maps produce  biconnected  graphs with two vertices only from biconnected  ones. Now, let $\vertex>1$. By assumption,
  the graph $\graph$ has at least one internal edge which
 does not belong to  any
 cycle.
   Therefore, it connects  two (distinct) vertices that
 must be connected together with only one internal edge. Let these vertices be  $v_i, v_j\in V$; $i,j\in\{1,\ldots,\vertex\}$ with $i\neq j$, for instance. Apply
 the $\Qione$-map to the graph $\graph$.
In particular, $\Qione(\graph)$ is a linear combination of graphs, each of which   is so that the   vertex $v_j$
 is connected with only one internal edge either to $v_i$ or $v_{\vertex+1}$
(but not to both). Clearly,  only
 cycles  containing the vertex $v_i$, are affected by the  $\Qione$-map.   That is, the vertex $v_j$
 cannot share a
 cycle with neither of  the vertices $v_i$ or $v_{\vertex+1}$. Hence,  the  graphs in $\Qione(\graph)$ are not biconnected. This completes the proof.
\end{proof}

\begin{lem}\label{lem:cyclestocycles}
Fix integers $\lp,\ext\ge0$ and $\vertex\ge 1$. Then, for all $i\in\{1,\ldots,\vertex\}$, the following statements hold:
\begin{enumerate}[(a)]\item $\Ti(\Q \Vb^{\vertex,\lp,\ext})\subseteq\Q\Vb^{\vertex,\lp+1,\ext}$;
\item $\Qcione(\Q\Vb^{\vertex,\lp,\ext})\subseteq\Q\Vb^{\vertex+1,\lp,\ext}$.\end{enumerate}
\end{lem}
\begin{proof}
(a) Clearly, the statement holds.  (b) Let $\graph=(V,K,E)$ 
together with the maps $\mapi$ and $\mape$, denote a graph in $\Vb^{\vertex, \lp,\ext}$.  Let $\mathcal{E}_i$ denote the set of ends of edges assigned to the vertex $v_i\in V$.
Apply the $\sci$-map to the graph $\graph$. 
 In particular, $\sci(\graph)$ is a linear combination of graphs, each of which  is produced from the graph $\graph$ by
transforming  one or more cycles containing the vertex $v_i$, into  paths
whose end point vertices are $v_i$ and $v_{\vertex+1}$. Moreover, every way to  assign the remaining ends of internal edges in the process, from $v_i$  to either $v_i$ or $v_{\vertex+1}$, defines new cycles.  Therefore, applying the $\Ein$-map to $\sci(\graph)$,
restores the broken cycles and yields biconnected graphs.  This completes the proof.
\end{proof}

%% file: loops.tex
\section{Arbitrary connected graphs}
\label{sec:loops}
The present section has substantial overlap with Section II  of  \cite{MeOe:loop}. Its main result is a  recursion formula to generate all connected graphs directly in the algebraic representation rooted in \cite{MeOe:npoint}. Here, we formulate
that formula on the level  of graphs. In a recursion step, the formula yields the linear combination of all graphs having the same vertex and cyclomatic numbers. Moreover,
 the sum of the 
 coefficients of all graphs in the same equivalence class, corresponds to  the inverse of the order of
their  group of automorphisms. Notice that  in \cite{MeOe:npoint,MeOe:loop}, the ordering of the vertices is not explicitly taken into account. That is, only one representative of each equivalence class is considered, the 
coefficient of such graph being  given by  the sum of the coefficients of all graphs in the same equivalence class.

\bigskip

We  use the $\Ti$ and $\Qione$-maps defined in the preceding section  to   recursively generate all connected  graphs.
\begin{thm}\label{thm:loops}
Fix an integer $\ext\ge0$.   
For all integers $\lp\ge0$ and $\vertex\ge 1$,   define $\Dl^{\vertex,\lp,\ext} \in \Q\Vc^{\vertex,\lp,\ext}$
by the following recursion relation:
\begin{itemize}
\item
$\Dl^{1,0,\ext}$ is a single vertex with $\ext$  external edges whose free ends are labeled $x_1,\dots,x_\ext$, and unit coefficient;
\item
\begin{equation}
\Dl^{\vertex,\lp,\ext}  \defeq
\frac{1}{\lp+\vertex-1}\left(\sum_{i=1}^{\vertex-1}\Qione(\Dl^{\vertex-1,\lp,\ext})+\frac{1}{\two}\sum_{i=1}^{\vertex}\Ti(\Dl^{\vertex,\lp-1,\ext})\right)\,.
 \label{eq:recloops}
\end{equation}
\end{itemize}
Then, for fixed values of  $\vertex$ and $\lp$, $\Dl^{\vertex,\lp,\ext}=\sum_{\graph\in \Vc^{\vertex,\lp,\ext}}\alpha_\graph\, \graph$; $\alpha_\graph\in\Q$ and $\alpha_\graph>0$; for all $\graph\in\Vc^{\vertex,\lp,\ext}$. Moreover,  $\sum_{\graph\in\C}\alpha_\graph=1/S^{\C}$, where $\C\subseteq\Vc^{\vertex,\lp,\ext}$ denotes an arbitrary equivalence 
class of graphs
and $S^{\C}$ denotes their symmetry factor.
\end{thm}
In the recursion equation above, the 
$\Ti$ summand does not appear when 
$\lp=0$.
In particular, for $\lp=0$,  formula (\ref{eq:recloops}) specializes to recursively generate all trees. Moreover, formula (\ref{eq:recloops}) is an instance of a double recursion. Therefore, its algorithmic implementation is that of any recursive function that makes two calls to itself, such as the defining recurrence of the binomial coefficients.

\begin{proof}
The proof is nearly the same to that of  Theorem 10 of \cite{MeOe:loop}. The procedure is also very analogous to the one given in \cite{MeOe:npoint}.  We translate every lemma  given in Section II of the former paper to the present setting.
\begin{lem}
\label{lem:allloops}
Fix integers $\ext,\lp\ge0$  and  $\vertex\ge1$. 
 Let $\Dl^{\vertex,\lp,\ext}=\sum_{\graph\in \Vc^{\vertex,\lp,\ext}}\alpha_\graph\, \graph \in \Q\Vc^{\vertex,\lp,\ext}$ be defined by formula (\ref{eq:recloops}). Then, $\alpha_\graph>0$ for all $\graph\in\Vc^{\vertex,\lp,\ext}$.
\end{lem}
\begin{proof}
The proof proceeds by induction on the internal edge  number $\edge$.  
Clearly, the statement holds for $\edge=0$.
We assume the result to hold for an arbitrary number of internal edges $\edge-1$. Let $\graph=(V,K,E)\in\Vc^{\vertex,\lp,\ext}$;  $E=\Int\cup\Ext$, $\edge=\mbox{card}(\Int)=\lp+\vertex-1$,
together with the maps $\mapi$ and $\mape$, denote a graph.   We show that the graph $\graph$  is generated by applying  the  $\Ti$-maps to graphs occurring in $\Dl^{\vertex,\lp-1,\ext}=\sum_{\graph^*\in\Vc^{\vertex,\lp-1,\ext}}\gamma_{\graph^*}\graph^*$; $\gamma_{\graph^*}\in\Q$, or the $\Qione$-maps   to graphs occurring in $\Dl^{\vertex-1,\lp,\ext}=\sum_{\graph'\in\Vc^{\vertex-1,\lp,\ext}}\beta_{\graph'}\graph'$; $\beta_{\graph'}\in\Q$:
\begin{enumerate}[(i)]\item  
 Suppose that the graph $\graph$ has at least one vertex with
one or more loops. Let this vertex be $v_i\in V$; $i\in\{1,\ldots,\vertex\}$, for instance. Erasing any loop,
yields a graph $\Ti^{-1}(\graph)\in\Q\Vc^{\vertex,\lp-1,\ext}$. By induction assumption, $\gamma_{\Ti^{-1}(\graph)}>0$. 
  Hence, applying the $\Ti$-map to  the graph $\Ti^{-1}(\graph)$ produces again the graph $\graph$. That is, $\alpha_{\graph}>0$. \item Suppose that the  graph $\graph$ has no loops.  There exists $i\in\{1,\ldots,\vertex-1\}$ so that $\{v_i,v_{\vertex}\}\in\mapi(\Int)$. Applying the $\cin$-map to the graph $\graph$ yields a graph $\cin(\graph)\in\Q\Vc^{\vertex-1,\lp,\ext}$. 
By induction assumption, $\beta_{\cin(\graph)}>0$.   Hence, applying the $\Qione$-map to the graph  $\cin(\graph)$,  produces a linear combination of graphs, one of which is the graph $\graph$.  That is, $\alpha_\graph>0$. \end{enumerate}  
\end{proof}
What remains in order to prove Theorem~\ref{thm:loops} is to show
that the sum of the coefficients of all graphs in the same equivalence class, is  given by the inverse of their symmetry factor. We start with a more
restricted result.
\begin{lem}
\label{lem:distloops} Fix integers $\lp\ge0$,  $\vertex\ge1$ and $\ext\ge \vertex$. 
 Let  $\C\subseteq \Vc^{\vertex,\lp,\ext}$ denote an equivalence class.  Let $\graph=(V,K,E)$; $E=\Int\cup\Ext$ 
together with the maps $\mapi$ and $\mape$, denote a graph in $\C$. Assume that $V\cap\mape(\Ext)=V$.  Let $\Dl^{\vertex,\lp,\ext}=\sum_{\graph\in \Vc^{\vertex,\lp,\ext}}\alpha_\graph\, \graph \in \Q\Vc^{\vertex,\lp,\ext}$ be defined by formula (\ref{eq:recloops}). Then, $\sum_{\graph\in\C}\alpha_{\graph}=1/S^{\C}$, where $S^{\C}$ denotes the symmetry factor of every graph in $\C$.
\end{lem}
\begin{proof}
The proof proceeds by induction on the internal edge  number $\edge$.  
Clearly, the
statement holds for $\edge=0$.  We assume the statement to  hold for a
general internal edge number $\edge-1$. Consider the graph $\graph=(V,K,E)\in\C$;   
$E=\Int\cup\Ext$, $\edge=\mbox{card}(\Int)=\lp+\vertex-1$, 
together with the maps $\mapi$ and $\mape$. 
 By Lemma~\ref{lem:allloops}, the coefficient of the graph $\graph$ in $\Dl^{\vertex,\lp,\ext}$ is positive, i.e., $\alpha_{\graph}>0$. 
 We proceed to show that
$\sum_{\graph\in\C}\alpha_{\graph}=1/S^{\C}$.  In particular, the graph $\graph\in\C$ is so that every one of its vertices has at least one (labeled) external edge. That is,  
the graph
$\graph$ has no non-trivial vertex
symmetries: $S^{\C}_{\text{vertex}}=1$. Hence,
$S^{\C}=S_{\text{edge}}^{\C}$ as any symmetry is an edge
symmetry.
We check from which graphs
with $\edge-1$ internal edges, the graphs in  the equivalence class   $\C$ are generated  by the recursion formula (\ref{eq:recloops}), and how many times they are generated. To this end, 
choose any one of the $\edge$ internal edges of the graph $\graph\in\C$:

\begin{enumerate}[(i)]
\item
If that internal edge is a loop, let this be assigned to the vertex $v_i\in V$; $i\in\{1,\ldots,\vertex\}$, for instance. Also, assume that the vertex $v_i$  has exactly $1\le\tau\le\lp$ loops as well as $x\ge1$ external edges whose free ends are labeled $x_{a_1},\dots,x_{a_x}$, with  $1\le a_1<\dots< a_x\le \ext$. 
 Erasing any one of these loops
yields a graph $\Ti^{-1}(\graph)$ whose
symmetry factor is related to that of $\graph\in\C$ via
$S^{\Ti^{-1}(\graph)}=S^{\C}/ (\two \tau)$. 
Let $\Dl^{\vertex,\lp-1,\ext}=\sum_{\graph^*\in\Vc^{\vertex,\lp-1,\ext}}\gamma_{\graph^*}\graph^*\,;\gamma_{\graph^*}\in\Q$.  
Also, let  $\A\subseteq\Vc^{\vertex,\lp-1,\ext}$ denote the  equivalence class containing $\Ti^{-1}(\graph)$.  The $\Ti$-map  produces the graph $\graph$ from the graph  $\Ti^{-1}(\graph)$  with coefficient $\alpha^*_{\graph}=
\gamma_{\Ti^{-1}(\graph)}\in\Q$. Each vertex of the graph $\Ti^{-1}(\graph)$ has at least one labeled external edge. Hence,  by induction assumption,  $\sum_{\graph^*\in\A}\gamma_{\graph^*}=1/S^{\Ti^{-1}(\graph)}=1/S^{\A}$. Now, take one graph (distinct from the graph $\graph$) in $\C$ in turn, choose one of the loops of the vertex having $x$ external edges whose free ends are labeled $x_{a_1},\dots,x_{a_x}$, and   repeat the procedure above. We obtain
\begin{eqnarray*}
\sum_{\graph\in\C}\alpha^*_{\graph}  =   \sum_{\graph^*\in\A}\gamma_{\graph^*}
 =  \frac{1}{S^{\A}}
 =   \frac{\two \tau}{S^{\C}} .
\end{eqnarray*}
  Therefore, the contribution to $\sum_{\graph\in\C}\alpha_{\graph}$ 
 is $\tau/ (\edge\cdot S^{\C})$. Distributing
this factor between the $\tau$ loops considered yields
$1/(\edge\cdot S^{\C})$ for each loop. 

\item 
If that internal edge is not a loop, let this be connected to the vertices, $v_i, v_j\in V$; $i,j\in\{1,\ldots,\vertex\}$ with $i<j$,   for instance. Also, assume  that $v_i$ has $\tau'\ge0$ loops as well as $r\ge1$ external edges whose free ends are labeled $x_{a_1},\dots,x_{a_r}$, with  $1\le a_1<\dots< a_r\le \ext$,  while $v_j$ has $\tau''\ge0$ loops as well as   $r'\ge1$ external edges whose  free ends are labeled
$x_{b_1},\dots,x_{b_{r'}}$, with  $1\le b_1<\dots< b_{r'}\le \ext$ and $a_z\neq b_{z'}$ for all $z\in \{1,\ldots,r\}, z'\in\{1,\ldots,r'\}$. Finally, assume that  the two vertices are connected together with $\multiedges\ge1$ internal edges, so that $1\le\tau'+\tau''+\multiedges\le\lp+1$.
Contracting any one of these internal edges,
yields a graph $\cij(\graph)$ whose $i$th
vertex, has $\mu\defeq \tau'+\tau''+\multiedges-1$ loops as well as $r+r'$   external edges whose free ends are labeled $x_{a_1},\dots,x_{a_r},x_{b_1},\dots,x_{b_{r'}}$. 
Consequently,
the symmetry factor  of the graph $\cij(\graph)$  is
related to that of $\graph\in\C$ via
$$
\frac{1}{2^\mu}\frac{1}{\mu!}S^{\cij(\graph) }
=\frac{1}{2^{\tau'}\tau'!}\frac{1}{2^{\tau''}\tau''!}\frac{1}{\multiedges!}
S^{\C}\,.$$
 Let $\Dl^{\vertex-1,\lp,\ext}=\sum_{\graph'\in\Vc^{\vertex-1,\lp,\ext}}\beta_{\graph'}\graph'\,;\beta_{\graph'}\in\Q$. 
 Let  $\B\subseteq\Vc^{\vertex-1,\lp,\ext}$ denote the  equivalence class containing $\cij(\graph)$. 
Applying the $\Qione$-map  to $\cij(\graph)$  yields a linear combination of graphs, one of which, 
 is isomorphic to the graph $\graph$. To calculate the coefficient $\alpha'_{\graph}\in\Q$ of  such graph 
in that linear combination, we need to count the number of different ways to distribute the $2\mu+r+r'$ ends of edges assigned to the  vertex $v_i$, 
between the two new ones, so that one vertex is assigned with $r$ external edges whose  free ends are labeled $x_{a_1},\dots,x_{a_r}$, as well as $\tau'$ loops, the other is assigned with $r'$ external edges whose  free ends  are labeled by $x_{b_1},\dots,x_{b_{r'}}$, as well as $\tau''$ loops, while the remaining $\multiedges-1$ internal edges are employed to connect the two vertices together. Now, 
there are two ways to assign the given $r$ external edges
to one vertex and  the given $r'$ external edges  to the other. Moreover, there are $\binom{\mu}{\tau'}=\frac{\mu!}{(\mu-\tau')!\tau'!}$ ways to assign both ends of $\tau'$ internal edges chosen among the $\mu$ internal edges in the process, to the vertex with the aforesaid $r$ external edges. Besides, there are  $\binom{\mu-\tau'}{\tau''}=\frac{(\mu-\tau'
)!}{(\mu-\tau-\tau'')!\tau''!}$ ways to assign both ends of $\tau''$ internal edges chosen among the $\mu-\tau'$ internal edges in the process, to the vertex with the aforesaid $r'$ external edges. Finally, 
there are two
ways to distribute one end of each of the remaining $\multiedges-1$ internal edges, per vertex. This yields $2^{\multiedges-1}$ ways to connect the two new vertices together with $\multiedges-1$ internal  edges. The final result is given by the product of all these factors. Hence, there are
\begin{equation*}
\label{eq:waysloops}
2\cdot 2^{\multiedges-1}\frac{\mu!}{(\mu-\tau')!\tau'!}\cdot\frac{(\mu-\tau')!}{(\mu-\tau'-\tau'')!\tau''!}=2^{\multiedges}\,\frac{\mu!}{\tau'!\,\tau''!\,(\multiedges-1)!}
\end{equation*}
ways to distribute the $2\mu+r+r'$ ends of edges between the two new vertices in order to produce a graph in the equivalence class $\C$. 
Hence, $\alpha'_{\graph}= 2^{\multiedges-1}\,\frac{\mu!}{\tau'!\,\tau''!\,(\multiedges-1)!}\beta_{\cij(\graph)}$. Each vertex of the graph $\cij(\graph)$  has at least one labeled external edge. Thus, by induction assumption,  $\sum_{\graph'\in\B}\beta_{\graph'}=1/S^{\cij(\graph)}=1/S^{\B}$. Now, take one graph (distinct from $\graph$) in $\C$ in turn, choose one of the internal edges connecting together the pair of vertices so that one vertex has $r$ external edges whose free ends are labeled $x_{a_1},\dots,x_{a_r}$, while the other has $r'$ external edges whose  free ends are labeled
$x_{b_1},\dots,x_{b_{r'}}$, and    repeat the procedure above. We obtain
\begin{eqnarray*}\nonumber
\sum_{\graph\in\C}\alpha'_{\graph} & = &  2^{\multiedges-1}\,\frac{\mu!}{\tau'!\,\tau''!\,(\multiedges-1)!}\sum_{\graph'\in\B}\beta_{\graph'}\\\nonumber
& = &  2^{\multiedges-1}\,\frac{\mu!}{\tau'!\,\tau''!\,(\multiedges-1)!}\frac{1}{S^{\B}}\\\nonumber
& = &  2^{\multiedges-1}\,\frac{\mu!}{\tau'!\,\tau''!\,(\multiedges-1)!}\cdot\frac{\tau'!\,\tau''!\,\multiedges!}{\mu!}\cdot\frac{1}{2^{\multiedges-1}S^{\C}}\\
& = & \frac{\multiedges}{S^{\C}}\,.\\
\end{eqnarray*}
  Therefore, the contribution to $\sum_{\graph\in\C}\alpha_{\graph}$ 
 is $\multiedges/ (\edge\cdot S^{\C})$.
Distributing
this factor between the $\multiedges$ internal edges considered yields
$1/(\edge\cdot S^{\C})$ for each edge.
\end{enumerate}
 
We conclude that every one of the $\edge$ internal edges of the graph $\graph$  contributes with a
factor of
$1/(\edge\cdot S^{\C})$ to  $\sum_{\graph\in\C}\alpha_{\graph}$. Hence, the overall contribution  is exactly $1/S^{\C}$. This
completes the proof.
\end{proof}

$\Dl^{\lp,\vertex,n}$ satisfies the following property:
\begin{lem}\label{lem:facloops}
Fix integers $\lp,\ext\ge 0$ and $\vertex\ge 1$. 
 Let $\Dl^{\vertex,\lp,\ext}=\sum_{\graph\in \Vc^{\vertex,\lp,\ext}}\alpha_\graph\, \graph \in \Q\Vc^{\vertex,\lp,\ext}$ be defined by formula (\ref{eq:recloops}). Let $K'\subset\mathcal{K}$ be a finite set so that $K\cap K'=\emptyset$.  Let $\Ext'\subseteq [K']^\two$; \emph{$\card(\Ext')=\ext'$}. Also, assume that  the elements of $\Ext'$ satisfy  $\{\legs_b,\legs_{b'}\}\cap\{\legs_{c},\legs_{c'}\}=\emptyset$.  Let  $L'=\{x_{\ext+1},\ldots,x_{\ext+\ext'}\}$ be a label set so that  $L\cap L'=\emptyset$. 
Let  $l':\Ext'\rightarrow [K',L']$ be a labeling of the elements of $\Ext'$. 
 Then,  $\Dl^{\vertex,\lp,\ext+\ext'}=\extb_{\Ext',V}(\Dl^{\vertex,\lp,\ext}).$
\end{lem}
\begin{proof}
Let $\Ext^*\defeq\Ext\cup\Ext'$ be the set of external edges  of all graphs occurring in   $\Dl^{\vertex-1,\lp,\ext+\ext'}\in\Q\Vc^{\vertex-1,\lp,\ext+\ext'}$.
Let $V'=\{v_1,\dots,v_{\vertex-1}\}$ be their vertex  set. 
Let   $\mathcal{E}_i$ be  the set of ends of edges assigned to the vertex  $v_i\in V'$;  $i\in\{1,\ldots, \vertex-1\}$. Let $\legsie$ be the subset of $\mathcal{E}_i$ whose elements are   ends of external edges. Also, let $\edgesie$ be the set of external edges assigned to the vertex $v_i$.  Finally, let $\edgesie^*\defeq\Ext'\cap\edgesie$; $\card(\edgesie^*)=\ext^*$, and $\legsie^*\defeq\legsie\cap\edgesie^*$. In this context, for all $i\in\{1,\ldots,\vertex-1\}$, the $\si$-maps yield  as $\si=\extb_{\edgesie^*,\{v_i,v_{\vertex}\}}\circ s'_{\mathcal{E}_i\backslash\legsie^*}$, where the maps $s'_{\mathcal{E}_i\backslash\legsie^*}:\Q\Vc^{\vertex-1,\lp,\ext+\ext'}\to\Q V^{\vertex,\lp-1,\ext+\ext'-\ext^*}\cup\Q V^{\vertex,\lp,\ext+\ext'-\ext^*}$  are required to produce graphs with external edge set $\Ext^*\backslash\edgesie^*$, from graphs with external edge set $\Ext^*$. Clearly, $\extb_{\Ext',V}=\extb_{\edgesie^*,\{v_i,v_\vertex\}}\circ\extb_{\Ext'\backslash\edgesie^*,V\backslash\{v_i,v_\vertex\}}$. Hence, the equality $\Dl^{\vertex,\lp,\ext+\ext'}=\extb_{\Ext',V}(\Dl^{\vertex,\lp,\ext})$ follows immediately from the recursive definition (\ref{eq:recloops}). 
 \end{proof}
We now proceed  to show that $\sum_{\graph\in\Vc^{\vertex,\lp,\ext}}\alpha_{\graph}=1/S^{\C}$, where $\C\subseteq\Vc^{\vertex,\lp,\ext}$ denotes any equivalence class.
\begin{lem}
\label{lem:anyclass} Fix integers $\lp\ge0$ and  $\vertex\ge1$. 
 Let  $\C\subseteq \Vc^{\vertex,\lp,\ext}$ denote an arbitrary equivalence class.
   Let $\Dl^{\vertex,\lp,\ext}=\sum_{\graph\in \Vc^{\vertex,\lp,\ext}}\alpha_\graph\, \graph \in \Q\Vc^{\vertex,\lp,\ext}$ be defined by formula (\ref{eq:recloops}). Then, $\sum_{\graph\in\C}\alpha_{\graph}=1/S^{\C}$, where $S^{\C}$ denotes the symmetry factor of every graph in $\C$.
\end{lem}
\begin{proof}
Choose a graph $\graph=(V,K,E)\in\C$; 
$E=\Int\cup\Ext$, 
together with the maps $\mapi$ and $\mape$. If $\mape(\Ext)=V$, 
 we
simply recall Lemma \ref{lem:distloops}. Thus, we may now assume that there exists  a set $V'\subseteq V$; $\mbox{card}(V')=\ext'$  so that $V'\cap\mape(\Ext)=\emptyset$.  Let $K'\subset\mathcal{K}$ be a finite set so that $K\cap K'=\emptyset$.  Also, let $\Ext'\subseteq [K']^\two$; $\card(\Ext')=\ext'$. Assume that the elements of $\Ext'$ satisfy  $\{\legs_b,\legs_{b'}\}\cap\{\legs_{c},\legs_{c'}\}=\emptyset$.  Also, let $L'=\{x_{\ext+1},\ldots,x_{\ext+\ext'}\}$ be a  label  set so that $L\cap L'=\emptyset$. Finally, let  $l':\Ext'\rightarrow [K',L']$ be a labeling of the elements of $\Ext'$. Now, apply an $\exta_{\Ext'}$-map to the graph $\graph$. Let $\D\subseteq\Vc^{\vertex,\lp,\ext+\ext'}$ denote the equivalence class containing $\exta_{\Ext'}(\graph)$. Let $\Dl^{\vertex,\lp,\ext+\ext'}=\sum_{\graph'\in\Vc^{\vertex,\lp,\ext+\ext'}}\beta_{\graph'}\graph'\,;\beta_{\graph'}\in\Q$.  
   By Lemma \ref{lem:distloops},  $\sum_{\graph'\in \D}\beta_{\graph'}= 1/S^{\exta_{\Ext'}(\graph)}= 1/S^{\D}$. 
Since, in general, the $\exta_{\Ext'}$-maps are not uniquely defined, assume that there are $T$ distinct maps $\exta^{(l)}_{\Ext'}$, $l\in\{1,\ldots,T\}$ so that  $\exta_{\Ext'}^{(l)}(\graph)\in\D$. Clearly, $\beta_{\exta^{(l)}_{\Ext'}(\graph)}=\alpha_{\graph}>0$. Therefore, by repeating the same procedure for every graph in $\C$ and recalling Lemma~\ref{lem:facloops}, we obtain 
\begin{eqnarray*}
\sum_{\graph'\in\D}\beta_{\graph'} = \sum_{l=1}^{T}\sum_{\graph\in\C} \beta_{\exta^{(l)}_{\Ext'}(\graph)}= T\sum_{\graph\in\C}\alpha_{\graph} =  \frac{1}{S^{\D}}\,.\\
\end{eqnarray*}
That is, $\sum_{\graph\in\C}\alpha_{\graph} =  1/(T\cdot S^{\D})$. Now, every map $\exta^{(l)}_{\Ext'}$ defines  a vertex
symmetry of the graph $\graph$.
This 
can have no more than
these vertex symmetries, since the vertices that already carry (labeled) external edges,
are distinguishable and thus held fixed under any symmetry. Hence, $S^{\graph}_{\text{vertex}}=S^{\C}_{\text{vertex}}=T$.
Moreover,
$S^{\D}=S^{\D}_{\text{edge}}=S^{\C}_{\text{edge}}$. Finally, from the identity
 $S^{\C}=S^{\C}_{\text{vertex}}\cdot
S^{\C}_{\text{edge}}$, follows that
$\sum_{\graph\in\C}\alpha_{\graph}=1/S^{\C}$. \end{proof} 
\end{proof}

\medskip

Appendix \ref{app:loops} shows the result of computing all mutually non isomorphic connected graphs without external edges as  contributions to $\Dl^{\vertex,\lp,0}$, for internal edge number $\edge=\lp+\vertex-1\le 3$.

%% file: biconn.tex
\section{Extensions}\label{sec:extensions}
We generalize the recursion formula (\ref{eq:recloops}) to   biconnected, simple and loopless connected graphs. These three results were not obtained in previous papers, using the Hopf algebraic approach given  in \cite{MeOe:npoint,MeOe:loop}.

\subsection{Biconnected graphs}
\label{sec:biconn}
By Lemmas \ref{lem:biconntobiconn} and \ref{lem:cyclestocycles},  Theorem \ref{thm:loops}  specializes to biconnected graphs by replacing the  $\Qione$-maps by  the $\Qcione$-maps in formula (\ref{eq:recloops}).
\begin{thm}\label{thm:biconn}
Fix an integer $\ext\ge0$.   
For all integers $\lp\ge0$ and $\vertex\ge 1$,   define $\Dp^{\vertex,\lp,\ext} \in \Q\Vb^{\vertex,\lp,\ext}$
by the following recursion relation:
\begin{itemize}
\item
$\Dp^{1,0,\ext}$ is a single vertex with $\ext$  external edges whose free ends are labeled $x_1,\dots,x_\ext$, and unit coefficient;
\item
\begin{equation}
\Dp^{\vertex,\lp,\ext}  \defeq
\frac{1}{\lp+\vertex-1}\left(\sum_{i=1}^{\vertex-1}\Qcione(\Dp^{\vertex-1,\lp,\ext})+\frac{1}{\two}\sum_{i=1}^\vertex\Ti(\Dp^{\vertex,\lp-1,\ext})\right)\,,\lp>0\,.
 \label{eq:recbiconn}
\end{equation}
\end{itemize}
Then, for fixed values of  $\vertex$ and $\lp$, $\Dp^{\vertex,\lp,\ext}=\sum_{\graph\in \Vb^{\vertex,\lp,\ext}}\alpha_\graph\, \graph$; $\alpha_\graph\in\Q$ and $\alpha_\graph>0;\,\,\forall\graph\in\Vb^{\vertex,\lp,\ext}$. Moreover,  $\sum_{\graph\in\C}\alpha_\graph=1/S^{\C}$, where $\C\subseteq\Vb^{\vertex,\lp,\ext}$ denotes an arbitrary equivalence 
class of graphs
and $S^{\C}$ denotes their symmetry factor.
\end{thm}
For $\lp=1$ and $\vertex>1$, formula (\ref{eq:recbiconn}) specializes to recursively generate a cycle  with $\vertex$ vertices, $\ext$ external edges whose free ends are labeled $x_1,\dots,x_\ext$,
and coefficient $1/{(2\vertex)}$, from a cycle 
with $\vertex-1$ vertices, the given external edges and coefficient $1/{(2(\vertex-1))}$.

\medskip

Appendix \ref{app:biconn} shows the result of computing all mutually non isomorphic biconnected graphs without external edges as  contributions to $\Dp^{\vertex,\lp,0}$, for internal edge number $\edge=\lp+\vertex-1\le 4$.

%% file: simple.tex
\subsection{Simple connected graphs}
\label{sec:simple}
We generalize  Theorem \ref{thm:loops} to  simple connected graphs. To this end, we combine the $\Qdione$-maps with the $\Ebij$-maps in formula (\ref{eq:recloops}).

\bigskip

\begin{thm}\label{thm:simple}
Fix an integer $\ext\ge0$.  
For all integers $\lp\ge0$ and $\vertex\ge 1$,   define  $\Dnew^{\vertex,\lp,\ext} \in \Q\Vs^{\vertex,\lp,\ext}$
by the following recursion relation:
\begin{itemize}
\item
$\Dnew^{1,0,\ext}$ is a single vertex with $\ext$  external edges whose free ends are labeled $x_1,\dots,x_\ext$, and unit coefficient;
\item
\begin{equation}
\Dnew^{\vertex,\lp,\ext}  \defeq
\frac{1}{\lp+\vertex-1}\left(\sum_{i=1}^{\vertex-1}\Qdione(\Dnew^{\vertex-1,\lp,\ext})+\sum_{i=1}^\vertex\sum_{j=1}^{i-1}\Ebij(\Dnew^{\vertex,\lp-1,\ext})\right)\,, \vertex>1\,.
 \label{eq:recsimple}
\end{equation}
\end{itemize}
Then, for fixed values of  $\vertex$ and $\lp$, $\Dnew^{\vertex,\lp,\ext}=\sum_{\graph\in \Vs^{\vertex,\lp,\ext}}\alpha_\graph\, \graph$; $\alpha_\graph\in\Q$ and $\alpha_\graph>0$; for all $\graph\in\Vs^{\vertex,\lp,\ext}$. Moreover,   $\sum_{\graph\in\C}\alpha_\graph=1/S^{\C}$, where $\C\subseteq\Vs^{\vertex,\lp,\ext}$ denotes an arbitrary equivalence 
class of graphs
and $S^{\C}$ denotes their symmetry factor.
\end{thm}
In the recursion equation above, the $\Ebij$ summand does not appear when $\lp=0$ and/or $\vertex=2$.

\begin{proof}
The proof is very analogous to that of Theorem \ref{thm:loops} . Actually, every lemma given in the preceding section remains valid  by replacing $\Dl^{\vertex, \lp,\ext}$ by $\Dnew^{\vertex, \lp,\ext}$. Here,  we only state and prove the two  lemmas corresponding to Lemmas \ref{lem:allloops} and  \ref{lem:distloops}.  The rest of the proof  is implied by analogy.
\begin{lem}
\label{lem:allsimple}
Fix integers $\ext,\lp\ge0$  and  $\vertex\ge1$. 
 Let $\Dnew^{\vertex,\lp,\ext}=\sum_{\graph\in \Vs^{\vertex,\lp,\ext}}\alpha_\graph\, \graph \in \Q\Vs^{\vertex,\lp,\ext}$ be defined by formula (\ref{eq:recsimple}). Then, $\alpha_\graph>0$ for all $\graph\in\Vs^{\vertex,\lp,\ext}$.
\end{lem}
\begin{proof}
The proof proceeds by induction on the internal edge  number $\edge$. Clearly, the statement holds for $\edge=0$.
We assume the result to hold for an arbitrary internal edge number $\edge-1$. Let $\graph=(V,K,E)\in\Vs^{\vertex,\lp,\ext}$; 
 $E=\Int\cup\Ext$, $\edge=\mbox{card}(\Int)=\lp+\vertex-1$, 
together with the maps $\mapi$ and $\mape$, denote a graph. We show that the graph $\graph$  is generated by applying  the  $\Ebij$-maps to graphs occurring in $\Dnew^{\vertex,\lp-1,\ext}=\sum_{\graph^*\in\Vs^{\vertex,\lp-1,\ext}}\gamma_{\graph^*}\graph^*$; $\gamma_{\graph^*}\in\Q$, or the $\Qdione$-maps   to graphs occurring in $\Dnew^{\vertex-1,\lp,\ext}=\sum_{\graph'\in\Vs^{\vertex-1,\lp,\ext}}\beta_{\graph'}\graph'$; $\beta_{\graph^*}\in\Q$:
\begin{enumerate}[(i)]\item   If $\lp=0$, by Lemma \ref{lem:allloops}, $\alpha_\graph>0$.
\item
If $\lp>0$, choose any one of the internal edges of the graph $\graph$ which belong at least to one cycle. Let this be connected to  the vertices $v_i, v_j$; $i, j\in\{1,\ldots,\vertex\}$ with $i\neq j$, for instance. 
 Erasing such internal edge,
yields a graph $\Eij^{-1}(\graph)\in\Q\Vs^{\vertex,\lp-1,\ext}$. By induction assumption, $\gamma_{\Eij^{-1}(\graph)}>0$. 
  Hence, applying the $\Ebij$-map to  $\Eij^{-1}(\graph)$ produces again the graph $\graph$. That is, $\alpha_{\graph}>0$.  \end{enumerate}
\end{proof}
\begin{lem}
\label{lem:distsimple}
 Fix integers $\lp\ge0$,  $\vertex\ge1$ and $\ext\ge \vertex$. 
 Let  $\C\subseteq \Vs^{\vertex,\lp,\ext}$ denote an equivalence class.  Let $\graph=(V,K,E)$;  
$E=\Int\cup\Ext$, 
together with the maps $\mapi$ and $\mape$, denote a graph in $\C$. Assume that  $V\cap\mape(\Ext)=V$.  Let $\Dnew^{\vertex,\lp,\ext}=\sum_{\graph\in \Vs^{\vertex,\lp,\ext}}\alpha_\graph\, \graph \in \Q\Vs^{\vertex,\lp,\ext}$ be defined by formula (\ref{eq:recsimple}). Then, $\sum_{\graph\in\C}\alpha_{\graph}=1$. 
\end{lem}
\begin{proof}The proof proceeds by induction on the internal edge  number: $\edge=\lp+\vertex-1$. Clearly, the
statement holds for $\edge=0$.  We assume the statement to  hold for a
general  internal edge number $\edge-1$. Consider the graph $\graph=(V,K,E)\in\C$;  
$E=\Int\cup\Ext$, $\edge=\mbox{card}(\Int)=\lp+\vertex-1$, 
together with the maps $\mapi$ and $\mape$. 
 By Lemma~\ref{lem:allsimple}, the coefficient of the graph $\graph\in\C$ is positive: $\alpha_{\graph}>0$. 
In particular, $S^{\C}=1$ as the graph $\graph$ is simple and every one of its vertices has at least one external edge. We proceed to show that
$\sum_{\graph\in\C}\alpha_{\graph}=1$. To this end,   
choose any one of the $\edge$ internal edges of the graph $\graph\in\C$:
\begin{enumerate}[(i)]
\item
If that internal edge does not belong to any cycle, by  Lemma \ref{lem:distloops}, it  contributes
with a factor of  $1/\edge$ to $\sum_{\graph\in\C}\alpha_{\graph}$. 

\item If that internal edge belongs at least to one cycle, let this be connected to the vertices $v_i,v_j\in V$; $i,j\in\{1,\ldots,\vertex\}$ with $i\neq j$, for instance. Also, assume   that $v_i$ has  $r\ge1$ external edges whose free ends are labeled $x_{a_1},\dots,x_{a_r}$, with  $1\le a_1<\dots< a_r\le \ext$,  while $v_j$ has    $r'\ge1$ external edges whose  free ends are labeled
$x_{b_1},\dots,x_{b_{r'}}$, with  $1\le b_1<\dots< b_{r'}\le \ext$ and $a_z\neq b_{z'}$ for all $z\in \{1,\ldots,r\}, z'\in\{1,\ldots,r'\}$. 
Erasing the given internal edge
yields a graph $\Eij^{-1}(\graph)$ so that 
$S^{\Eij^{-1}(\graph)}=1$. 
Let $\Dnew^{\vertex,\lp-1,\ext}=\sum_{\graph^*\in\Vs^{\vertex,\lp-1,\ext}}\gamma_{\graph^*}\graph^*\,;\gamma_{\graph^*}\in\Q$. 
Also, let  $\A\subseteq\Vs^{\vertex,\lp-1,\ext}$ denote the equivalence class containing $\Eij^{-1}(\graph)$.  The $\Ebij$-map  produces the graph $\graph$ from the graph  $\Eij^{-1}(\graph)$  with coefficient $\alpha^*_{\graph}=
\gamma_{\Eij^{-1}(\graph)}\in\Q$. Each vertex of the graph $\Eij^{-1}(\graph)$ has at least one labeled external edge. Hence,  by induction assumption,  $\sum_{\graph^*\in\A}\gamma_{\graph^*}=1$. Now, take one graph (distinct from $\graph$) in $\C$ in turn, choose the internal edge connected to the  vertex having $r$ external edges whose free ends are labeled $x_{a_1},\dots,x_{a_r}$, and to the vertex having $r'$ external edges whose  free ends are labeled
$x_{b_1},\dots,x_{b_{r'}}$, and    repeat the procedure above. We obtain
\begin{eqnarray*}
\sum_{\graph\in\C}\alpha^*_{\graph}  =   \sum_{\graph^*\in\A}\gamma_{\graph^*}
 =   1 .
\end{eqnarray*}
  Therefore, the contribution to $\sum_{\graph\in\C}\alpha_{\graph}$ 
 is $1/\edge$. 
\end{enumerate}

\smallskip

 We conclude that every one of the $\edge$ internal edges of the graph $\graph$  contributes with a
factor of
$1/\edge$ to  $\sum_{\graph\in\C}\alpha_{\graph}$. Hence, the overall contribution  is exactly $1$. This
completes the proof.
\end{proof}
\end{proof}

\medskip

Appendix \ref{app:simple} shows the result of computing all mutually non isomorphic simple connected graphs without external edges as  contributions to $\Dnew^{\vertex,\lp,0}$, for internal edge number $\edge=\lp+\vertex-1\le 5$.

%% file: noloops.tex
\subsection{Loopless connected
graphs}
\label{sec:noloops}
The present section presents two algorithms to generate all loopless connected graphs. The second one is amenable for direct implementation via Hopf algebras in the sense of \cite{MeOe:npoint,MeOe:loop}.  
\subsubsection{Main recursion formula}
We generalize  Theorem \ref{thm:loops} to loopless connected graphs. To this end, we replace the $\Ti$-maps by the $\Eaij$-maps in formula (\ref{eq:recloops}).

\bigskip

\begin{thm}\label{thm:noloops}
Fix an integer $\ext\ge0$.   
For all integers $\lp\ge0$ and $\vertex\ge 1$,   define  $\Ds^{\vertex,\lp,\ext} \in \Q\Vm^{\vertex,\lp,\ext}$ by the following recursion relation:
\begin{itemize}
\item
$\Ds^{1,0,\ext}$ is a single vertex with $\ext$  external edges whose free ends are labeled $x_1,\dots,x_\ext$, and unit coefficient;
\item
\begin{equation}
\Ds^{\vertex,\lp,\ext}  \defeq
\frac{1}{\lp+\vertex-1}\left(\sum_{i=1}^{\vertex-1}\Qione(\Ds^{\vertex-1,\lp,\ext})+\sum_{i=1}^\vertex\sum_{j=1}^{i-1}\Eaij(\Ds^{\vertex,\lp-1,\ext})\right)\,, \vertex>1\,.
 \label{eq:recnoloops}
\end{equation}
\end{itemize}
Then, for fixed values of  $\vertex$ and $\lp$, $\Ds^{\vertex,\lp,\ext}=\sum_{\graph\in \Vm^{\vertex,\lp,\ext}}\alpha_\graph\, \graph$; $\alpha_\graph\in\Q$ and $\alpha_\graph>0$ for all $\graph\in\Vm^{\vertex,\lp,\ext}$. Moreover,   $\sum_{\graph\in\C}\alpha_\graph=1/S^{\C}$, where $\C\subseteq\Vm^{\vertex,\lp,\ext}$ denotes an arbitrary equivalence
class of graphs
and $S^{\C}$ denotes their symmetry factor.
\end{thm}
In the recursion equation above, the  $\Eaij$ summand does not appear  when $\lp=0$.

\begin{proof}
As in the preceding section, every lemma given in Section \ref{sec:loops} holds when stated for $\Ds^{\vertex,\lp,\ext}$. Hence, we restrict the proof of Theorem \ref{thm:noloops} to  the following two lemmas.
\begin{lem}
\label{lem:allnoloops}
Fix integers $\ext,\lp\ge0$  and  $\vertex\ge1$. 
 Let $\Ds^{\vertex,\lp,\ext}=\sum_{\graph\in \Vm^{\vertex,\lp,\ext}}\alpha_\graph\, \graph \in \Q\Vm^{\vertex,\lp,\ext}$ be defined by formula (\ref{eq:recnoloops}). Then, $\alpha_\graph>0$ for all $\graph\in\Vm^{\vertex,\lp,\ext}$.
\end{lem}
\begin{proof}The proof proceeds by induction on the internal edge  number $\edge$. Clearly, the statement holds for $\edge=0$.
We assume the result to hold for an arbitrary  internal edge number $\edge-1$. Let $\graph=(V,K,E)\in\Vm^{\vertex,\lp,\ext}$; 
$E=\Int\cup\Ext$, $\edge=\mbox{card}(\Int)=\lp+\vertex-1$,
together with the maps $\mapi$ and $\mape$, denote a graph. We show that the graph $\graph$  is generated by applying  the  $\Eaij$-maps to graphs occurring in $\Ds^{\vertex,\lp-1,\ext}=\sum_{\graph^*\in\Vm^{\vertex,\lp-1,\ext}}\gamma_{\graph^*}\graph^*$; $\gamma_{\graph^*}\in\Q$, or the $\Qione$-maps   to graphs occurring in $\Ds^{\vertex-1,\lp,\ext}=\sum_{\graph'\in\Vm^{\vertex-1,\lp,\ext}}\beta_{\graph'}\graph'$; $\beta_{\graph^*}\in\Q$:
\begin{enumerate}[(i)]\item   Suppose that the graph $\graph$ has no multiple edges. By Lemma \ref{lem:allloops}, $\alpha_\graph>0$. 
\item Suppose that the graph $\graph$ has at least one pair of vertices, say, $v_i,v_j\in V$; $i,j\in\{1,\ldots,\vertex\}$ with $i\neq j$, connected together by multiple edges.
 Erasing any one of those edges,
yields a graph $\Eij^{-1}(\graph)\in\Q\Vm^{\vertex,\lp-1,\ext}$. By induction assumption, $\gamma_{\Eij^{-1}(\graph)}>0$.   Hence, applying the $\Eaij$-map to  $\Eij^{-1}(\graph)$ produces again the graph $\graph$. That is, $\alpha_{\graph}>0$.
 \end{enumerate}
\end{proof}
\begin{lem}
\label{lem:distnoloops}
Fix integers $\lp\ge0$,  $\vertex\ge1$ and $\ext\ge \vertex$. 
 Let  $\C\subseteq \Vm^{\vertex,\lp,\ext}$ denote an equivalence class.  Let $\graph=(V,K,E)$;  
$E=\Int\cup\Ext$, 
together with the maps $\mapi$ and $\mape$, denote a graph in $\C$. Assume that $V\cap\mape(\Ext)=V$.  Let $\Ds^{\vertex,\lp,\ext}=\sum_{\graph\in \Vm^{\vertex,\lp,\ext}}\alpha_\graph\, \graph \in \Q\Vm^{\vertex,\lp,\ext}$ be defined by formula (\ref{eq:recnoloops}). Then, $\sum_{\graph\in\C}\alpha_{\graph}=1/S^{\C}$, where $S^{\C}$ denotes the symmetry factor of every graph in $\C$.
\end{lem}
\begin{proof}The proof proceeds by induction on the internal edge  number $\edge$. Clearly, the
statement holds for $\edge=0$.  We assume the statement to  hold for a
general number of internal edges $\edge-1$. Consider the graph $\graph=(V,K,E)\in\C$; 
$E=\Int\cup\Ext$, $\edge=\mbox{card}(\Int)=\lp+\vertex-1$,
together with the maps $\mapi$ and $\mape$. 
 By Lemma~\ref{lem:allnoloops}, the coefficient of the graph  $\graph\in\C$ is positive: $\alpha_{\graph}>0$.
 We proceed to show that
$\sum_{\graph\in\C}\alpha_{\graph}=1/S^{\C}$. To this end, 
choose any one of the $\edge$ internal edges of the graph $\graph\in\C$:
\begin{enumerate}[(i)]
\item
If that internal edge is the only one connecting a given pair of vertices together,  by  Lemma \ref{lem:distloops}, it  contributes
with a factor of  $1/(\edge\cdot S^{\C})$ to $\sum_{\graph\in\C}\alpha_{\graph}$. 

\item If that internal edge is one of the, say, $1<\rho\le\lp+1$, multiple edges connecting together the vertices  $v_i,v_j\in V$; $i, j\in\{1,\ldots,\vertex\}$ with $i\neq j$, for instance,  assume  that $v_i$ has  $r\ge1$ external edges whose free ends are labeled $x_{a_1},\dots,x_{a_r}$, with  $1\le a_1<\dots< a_r\le \ext$,  while $v_j$ has    $r'\ge1$ external edges whose  free ends are labeled
$x_{b_1},\dots,x_{b_{r'}}$, with  $1\le b_1<\dots< b_{r'}\le \ext$ and $a_z\neq b_{z'}$ for all $z\in \{1,\ldots,r\}, z'\in\{1,\ldots,r'\}$.
Erasing any one of the given internal edges
yields a graph $\Eij^{-1}(\graph)$ whose
symmetry factor is related to that of the graph $\graph\in\C$ via
$S^{\Eij^{-1}(\graph)}=S^{\C}/ \multiedges$.
Let $\Ds^{\vertex,\lp-1,\ext}=\sum_{\graph^*\in\Vm^{\vertex,\lp-1,\ext}}\gamma_{\graph^*}\graph^*\,;\gamma_{\graph^*}\in\Q$. 
Also, let  $\A\subseteq\Vm^{\vertex,\lp-1,\ext}$ denote the equivalence class containing $\Eij^{-1}(\graph)$.  The $\Eaij$-map  produces the graph $\graph$ from the graph  $\Eij^{-1}(\graph)$  with coefficient $\alpha^*_{\graph}=
\gamma_{\Eij^{-1}(\graph)}\in\Q$. Each vertex of the graph $\Eij^{-1}(\graph)$ has at least one labeled external edge. Hence,  by induction assumption,  $\sum_{\graph^*\in\A}\gamma_{\graph^*}=1/S^{\Eij^{-1}(\graph)}=1/S^{\A}$. Now, take one graph (distinct from $\graph$) in $\C$ in turn, choose one of the internal edges connecting together the pair of vertices so that one vertex has $r$ external edges whose free ends are labeled $x_{a_1},\dots,x_{a_r}$, while the other has $r'$ external edges whose  free ends are labeled
$x_{b_1},\dots,x_{b_{r'}}$, and    repeat the procedure above. We obtain
\begin{eqnarray*}
\sum_{\graph\in\C}\alpha^*_{\graph}  =   \sum_{\graph\in\A}\gamma_{\graph^*}
 =   \frac{1}{S^{\A}}
 =   \frac{\rho}{S^{\C}}\,.
\end{eqnarray*}
  Therefore, the contribution to $\sum_{\graph\in\C}\alpha_{\graph}$
 is $\rho/ (\edge\cdot S^{\C})$. Distributing
this factor between the $\rho$ internal edges considered yields
$1/(\edge\cdot S^{\C})$ for each internal edge. 
\end{enumerate}

\smallskip

We conclude that every one of the $\edge$ internal edges of the graph $\graph$  contributes with a
factor of
$1/(\edge\cdot S^{\C})$ to  $\sum_{\graph\in\C}\alpha_{\graph}$. Hence, the overall contribution  is exactly $1/S^{\C}$. This
completes the proof.
\end{proof}
\end{proof}

%% file: theta.tex
\subsubsection{Alternative recursion formula}
\label{sec:theta}
We present an alternative  recursion formula for loopless connected graphs. The underlying algorithm is amenable to direct implementation using the algebraic representation of graphs given in \cite{MeOe:npoint, MeOe:loop}.

\bigskip

\begin{thm}\label{thm:theta}
Fix an integer $\ext\ge0$.   
For all integers $\lp\ge0$ and $\vertex\ge 1$,   define  $\Dsnew^{\vertex,\lp,\ext} \in \Q\Vm^{\vertex,\lp,\ext}$
by the following recursion relation:
\begin{itemize}
\item
$\Dsnew^{1,0,\ext}$ is a single vertex with $\ext$  external edges whose free ends are labeled $x_1,\dots,x_\ext$, and unit coefficient;
\item
$\Dsnew^{1,\lp,\ext}  \defeq  0$, $\lp>0$;
\item
\begin{equation}
\Dsnew^{\vertex,\lp,\ext}  \defeq \frac{1}{\lp+\vertex-1}\sum_{
\multiedges=1}^{\lp+1}
\sum_{i=1}^{\vertex-1}\Qi(\Dsnew^{\vertex-1,\lp+1-\multiedges,\ext}), \vertex>1\,.\label{eq:rectheta}
\end{equation}
\end{itemize}
Then, for fixed values of  $\vertex$ and $\lp$, $\Dsnew^{\vertex,\lp,\ext}=\sum_{\graph\in \Vm^{\vertex,\lp,\ext}}\alpha_\graph\, \graph$; $\alpha_\graph\in\Q$ and $\alpha_\graph>0$; for all $\graph\in\Vm^{\vertex,\lp,\ext}$. Moreover,   $\sum_{\graph\in\C}\alpha_\graph=1/S^{\C}$, where $\C\subseteq\Vm^{\vertex,\lp,\ext}$ denotes an arbitrary equivalence of graphs 
class of graphs
and $S^{\C}$  denotes their symmetry factor.
\end{thm}
\begin{proof}
As in the previous sections, every lemma given in Section \ref{sec:loops} holds for $\Dsnew^{\vertex, \lp,\ext}$.  Here, details are only given for the following two lemmas. 

\begin{lem}
\label{lem:alltheta}
Fix integers $\ext,\lp\ge0$  and  $\vertex\ge1$. 
 Let $\Dsnew^{\vertex,\lp,\ext}=\sum_{\graph\in \Vm^{\vertex,\lp,\ext}}\alpha_\graph\, \graph \in \Q\Vm^{\vertex,\lp,\ext}$ be defined by formula (\ref{eq:rectheta}). Then, $\alpha_\graph>0$ for all $\graph\in\Vm^{\vertex,\lp,\ext}$.
\end{lem}
\begin{proof} 
The proof proceeds by induction on the internal edge  number $\edge$. Clearly, the
statement holds for $\edge=0$.  
We assume the statement to hold for any internal edge number
smaller than a fixed $\edge\ge 1$. Let $\graph=(V,K,E)\in\Vm^{\vertex,\lp,\ext}$; 
$E=\Int\cup\Ext$, $\edge=\mbox{card}(\Int)=\lp+\vertex-1$,
together with the maps $\mapi$ and $\mape$, denote a graph.   We show that the graph $\graph$  is generated by applying   the $\Qi$-maps   to graphs occurring in $\Dsnew^{\vertex-1,\lp+1-\multiedges,\ext}=\sum_{\graph'\in\Vm^{\vertex-1,\lp+1-\multiedges,\ext}}\beta_{\graph'}\graph'$; $\beta_{\graph'}\in\Q$:  There exists $i\in\{1,\ldots,\vertex-1\}$ so that $\{v_i,v_{\vertex}\}\in\mapi(\Int)$. Assume that the two vertices are connected together with $1\le\multiedges\le\lp+1$ internal edges.  Applying the map $\cin\circ\Einn^{1-\multiedges}$ to the graph $\graph$ yields a graph $\graph'\defeq(\cin\circ\Einn^{1-\multiedges})(\graph)\in\Q\Vm^{\vertex-1,\lp-\multiedges+1,\ext}$. 
By induction assumption, $\beta_{\graph'}>0$.   Hence, applying the $\Qi$-map to  the graph $\graph'$,  produces a linear combination of graphs, one of which is the graph $\graph$.  That is, $\alpha_\graph>0$. This completes the proof.
\end{proof}
\begin{lem}
\label{lem:disttheta}
Fix integers $\lp\ge0$,  $\vertex\ge1$ and $\ext\ge \vertex$. 
 Let  $\C\subseteq \Vm^{\vertex,\lp,\ext}$ denote an equivalence class.  Let $\graph=(V,K,E)$;  
$E=\Int\cup\Ext$, 
together with the maps $\mapi$ and $\mape$, denote a graph in $\C$. Assume that $V\cap\mape(\Ext)=V$.  Let $\Dsnew^{\vertex,\lp,\ext}=\sum_{\graph\in \Vm^{\vertex,\lp,\ext}}\alpha_\graph\, \graph \in \Q\Vm^{\vertex,\lp,\ext}$ be defined by formula (\ref{eq:rectheta}). Then, $\sum_{\graph\in\C}\alpha_{\graph}=1/S^{\C}$, where $S^{\C}$ denotes the symmetry factor of every graph in $\C$.
\end{lem}
\begin{proof}The proof proceeds by induction on the internal edge  number $\edge$. Clearly, the
statement holds for $\edge=0$.  
We assume the statement to hold for any internal edge number
smaller than a fixed $\edge\ge 1$. Consider the graph $\graph=(V,K,E)\in\C$;  
$E=\Int\cup\Ext$; $\edge=\mbox{card}(\Int)=\lp+\vertex-1$, 
together with the maps $\mapi$ and $\mape$. 
 By Lemma \ref{lem:alltheta}, the coefficient of the graph $\graph\in\C$ is positive: $\alpha_{\graph}>0$. 
 We proceed to show that
$\sum_{\graph\in\C}\alpha_{\graph}=1/S^{\C}$: 
Choose any one of the $\edge$ internal edges of the graph $\graph\in\C$.  Let this be  connected to the vertices $v_i,v_j\in V$; $i, j\in\{1,\ldots,\vertex\}$ with $i< j$, for instance. Also, assume  that $v_i$ has  $r\ge1$ external edges whose free ends are labeled $x_{a_1},\dots,x_{a_r}$, with  $1\le a_1<\dots< a_r\le \ext$,  while $v_j$ has    $r'\ge1$ external edges whose  free ends are labeled
$x_{b_1},\dots,x_{b_{r'}}$, with  $1\le b_1<\dots< b_{r'}\le \ext$ and $a_z\neq b_{z'}$ for all $z\in \{1,\ldots,r\}, z'\in\{1,\ldots,r'\}$. Finally, assume that there are $1\le\rho\le\lp+1$ multiple edges connecting the vertices $v_i$ and $v_j$ together.
Erasing $\multiedges-1$ internal edges and contracting the final one,
yields a graph $\graph'=(\cij\circ\Eij^{1-\multiedges})(\graph)$ whose symmetry factor  is
related to that of the graph $\graph\in\C$ via
$
S^{\graph' }
=\frac{1}{
\multiedges!}
S^{\C}\,.$
 Let $\Dsnew^{\vertex-1,\lp+1-\multiedges,\ext}=\sum_{\graph'\in\Vm^{\vertex-1,\lp+1-\multiedges,\ext}}\beta_{\graph'}\graph'\,;\beta_{\graph'}\in\Q$. 
 Let  $\B\subseteq\Vm^{\vertex-1,\lp+1-\multiedges,\ext}$ denote the equivalence class containing $\graph'$. 
Applying the $\Qi$-map  to the graph $\graph'$ produces a graph isomorphic to the graph $\graph$ with coefficient $\alpha'_\graph=\frac{\beta_{\graph'}}{(\multiedges-1)!}\in\Q$. Each vertex of the graph $\graph'$  has at least one labeled external edge. Hence, by induction assumption,  $\sum_{\graph'\in\B}\beta_{\graph'}=1/S^{\graph'}=1/S^{\B}$. Now, take one graph (distinct from $\graph$) in $\C$ in turn, choose one of the internal edges connected to the vertex having $r$ external edges whose free ends are labeled $x_{a_1},\dots,x_{a_r}$, and to the vertex having $r'$ external edges whose  free ends are labeled
$x_{b_1},\dots,x_{b_{r'}}$, and    repeat the procedure above. We obtain
\begin{eqnarray*}
\sum_{\graph\in\C}\alpha'_{\graph}  =   \frac{1}{(\multiedges-1)!}\sum_{\graph'\in\B}\beta_{\graph'}
 =   \frac{1}{(\multiedges-1)!S^{\B}}
 =  \frac{\multiedges}{S^{\C}} .
\end{eqnarray*}
  Therefore, the contribution to $\sum_{\graph\in\C}\alpha_{\graph}$ 
 is $\multiedges/ (\edge\cdot S^{\C})$
Distributing
this factor between the $\multiedges$ internal edges considered yields
$1/(\edge\cdot S^\C)$ for each edge. 
\smallskip

We conclude that every one of the $\edge$ internal edges of the graph $\graph$  contributes with a
factor of
$1/(\edge\cdot S^{\C})$ to  $\sum_{\graph\in\C}\alpha_{\graph}$. Hence, the overall contribution  is exactly $1/S^{\C}$. This
completes the proof.
\end{proof}
\end{proof}

\medskip

Appendix \ref{app:noloops} shows the result of computing all mutually non isomorphic connected graphs without external edges as  contributions to $\Ds^{\vertex,\lp,0}$ or $\Dsnew^{\vertex,\lp,0}$, for internal edge number $\edge=\lp+\vertex-1\le 4$.

%% file: algo.tex
\subsubsection{Algorithmic considerations}
\label{sec:algo}
The results of the present section can be seen as an extension of those
of Section IV of \cite{MeOe:loop}, to loopless connected graphs.

\bigskip

The two algorithms  underlying the recursive definitions
(\ref{eq:recloops}) and (\ref{eq:rectheta}) given in Sections \ref{sec:loops} and \ref{sec:theta}, respectively,    are amenable for direct implementation using the Hopf algebraic representation of graphs given in \cite{MeOe:npoint,MeOe:loop}. This representation
can be used directly and efficiently in implementing concrete
calculations of  graphs.  

 An important aspect for the
efficiency of concrete calculations is to discard graphs that do not
contribute. For instance, assume that one is only interested in calculating loopless graphs so that all vertices have a minimum degree, say, $\nu\ge \lp+1$.  In particular, 
in the recursive definition (\ref{eq:rectheta}), the number of ends of edges assigned to a
vertex changes after applying the $\Qi$-maps. The only
graphs with degree $1\le\multiedges\le\lp+1$ vertices are those produced by 
the $\si$-maps when one of the new vertices receives no ends of edges at all. This,
 thus, acquires degree $\multiedges$ after being  connected  to the other vertex  with $\multiedges$ internal edges.
 Hence,  to eliminate the irrelevant graphs with degree $\nu'< \nu$
vertices, replace the  $\Qi$-maps by ${\Qi}_{\ge \nu}\defeq \frac{1}{\two(\multiedges-1)!}\Ein^{\multiedges}\circ\si^{(\multiedges)}$ in formula (\ref{eq:rectheta}), where the $\si^{(\multiedges)}$-maps are required to partition the set of ends of edges assigned to the vertex $v_i$, $\mathcal{E}_i$, into two sets whose cardinality is equal or greater than $\nu-\multiedges$.

When considering loop graphs as
well,
 we can no longer
globally restrict the image of the maps $\Qione=\Ein\circ\si$ in formula (\ref{eq:recloops}). However, if we are interested in
graphs only up to a maximal  cyclomatic number, say, $\lp'$, we may still restrict the
partitions of the set $\mathcal{E}_i$ as part of the definition of the  $\Qione$-maps, in certain instances. These are precisely the
instances when a later application of  the $\Ti$-maps to a graph cannot occur,
i.e., when the graph has already the maximal cyclomatic number $\lp'$.

%% file: acknowledge.tex
\subsection*{Acknowledgments}
The author
is  very grateful to Edita Pelantov\'a for several corrections and suggestions  to Sections \ref{sec:basics}, \ref{sec:linear} and \ref{sec:loops}  of
the present manuscript, and also for enlightening discussions on the standard definitions and terminology of graph theory.
The author  
would also like to thank  Robert Oeckl for important corrections and suggestions to  the first version of
the  manuscript, and  Brigitte Hiller for helpful remarks.  
The research was supported in part by the
Czech Ministry of Education, Youth
and Sports within the project LC06002.

%% file: apploops.tex
\section{}\label{app:loops}
This appendix is the same as that of \cite{MeOe:loop}. It shows  $\Dl^{\vertex,\lp,0}$ up to order $\vertex+\lp\le4$  and computed via formula  (\ref{eq:recloops}).   All graphs in the same equivalence class are identified as the same.   The coefficients in front  of graphs are the
inverses of the orders of their
groups of automorphisms.
\vspace{1cm}\\
\input{figures/n1k0.eepic}
\vspace{1cm}\\
\input{figures/n2k0.eepic}
\vspace{1cm}\\
\input{figures/n1k1.eepic}
\vspace{1cm}\\
\input{figures/n3k0.eepic}
\vspace{1cm}\\
\input{figures/n2k1.eepic}
\vspace{1cm}\\
\input{figures/n1k2.eepic}
\vspace{1cm}\\
\input{figures/n4k0.eepic}
\vspace{1cm}\\
\input{figures/n3k1.eepic}
\vspace{1cm}\\
\input{figures/n2k2.eepic}
\vspace{1cm}\\
\input{figures/n1k3.eepic}

%% file: figures/n1k0.eepic
\setlength{\unitlength}{0.00083333in}
\begingroup\makeatletter\ifx\SetFigFontNFSS\undefined%
\gdef\SetFigFontNFSS#1#2#3#4#5{%
  \reset@font\fontsize{#1}{#2pt}%
  \fontfamily{#3}\fontseries{#4}\fontshape{#5}%
  \selectfont}%
\fi\endgroup%
{\renewcommand{\dashlinestretch}{30}
\begin{picture}(1510,276)(0,-10)
\texture{44555555 55aaaaaa aa555555 55aaaaaa aa555555 55aaaaaa aa555555 55aaaaaa 
	aa555555 55aaaaaa aa555555 55aaaaaa aa555555 55aaaaaa aa555555 55aaaaaa 
	aa555555 55aaaaaa aa555555 55aaaaaa aa555555 55aaaaaa aa555555 55aaaaaa 
	aa555555 55aaaaaa aa555555 55aaaaaa aa555555 55aaaaaa aa555555 55aaaaaa }
\put(1427,142){\shade\ellipse{150}{150}}
\put(1427,142){\ellipse{150}{150}}
\put(15,78){\makebox(0,0)[lb]{\smash{{\SetFigFontNFSS{12}{14.4}{\rmdefault}{\mddefault}{\updefault}$\vertex=1,\lp=0$}}}}
\end{picture}
}

%% file: figures/n2k0.eepic
\setlength{\unitlength}{0.00083333in}
\begingroup\makeatletter\ifx\SetFigFontNFSS\undefined%
\gdef\SetFigFontNFSS#1#2#3#4#5{%
  \reset@font\fontsize{#1}{#2pt}%
  \fontfamily{#3}\fontseries{#4}\fontshape{#5}%
  \selectfont}%
\fi\endgroup%
{\renewcommand{\dashlinestretch}{30}
\begin{picture}(2119,281)(0,-10)
\texture{44555555 55aaaaaa aa555555 55aaaaaa aa555555 55aaaaaa aa555555 55aaaaaa 
	aa555555 55aaaaaa aa555555 55aaaaaa aa555555 55aaaaaa aa555555 55aaaaaa 
	aa555555 55aaaaaa aa555555 55aaaaaa aa555555 55aaaaaa aa555555 55aaaaaa 
	aa555555 55aaaaaa aa555555 55aaaaaa aa555555 55aaaaaa aa555555 55aaaaaa }
\path(1749,141)(2049,141)
\path(1749,141)(2049,141)
\put(2036,142){\shade\ellipse{150}{150}}
\put(2036,142){\ellipse{150}{150}}
\put(1760,146){\shade\ellipse{150}{150}}
\put(1760,146){\ellipse{150}{150}}
\put(15,78){\makebox(0,0)[lb]{\smash{{\SetFigFontNFSS{12}{14.4}{\rmdefault}{\mddefault}{\updefault}$\vertex=2\,,\lp=0$}}}}
\put(1420,83){\makebox(0,0)[lb]{\smash{{\SetFigFontNFSS{12}{14.4}{\rmdefault}{\mddefault}{\updefault}$\frac{1}{2}$}}}}
\end{picture}
}

%% file: figures/n1k1.eepic
\setlength{\unitlength}{0.00083333in}
\begingroup\makeatletter\ifx\SetFigFontNFSS\undefined%
\gdef\SetFigFontNFSS#1#2#3#4#5{%
  \reset@font\fontsize{#1}{#2pt}%
  \fontfamily{#3}\fontseries{#4}\fontshape{#5}%
  \selectfont}%
\fi\endgroup%
{\renewcommand{\dashlinestretch}{30}
\begin{picture}(1868,394)(0,-10)
\put(1785,250){\ellipse{120}{246}}
\texture{44555555 55aaaaaa aa555555 55aaaaaa aa555555 55aaaaaa aa555555 55aaaaaa 
	aa555555 55aaaaaa aa555555 55aaaaaa aa555555 55aaaaaa aa555555 55aaaaaa 
	aa555555 55aaaaaa aa555555 55aaaaaa aa555555 55aaaaaa aa555555 55aaaaaa 
	aa555555 55aaaaaa aa555555 55aaaaaa aa555555 55aaaaaa aa555555 55aaaaaa }
\put(1785,144){\shade\ellipse{150}{150}}
\put(1785,144){\ellipse{150}{150}}
\put(15,78){\makebox(0,0)[lb]{\smash{{\SetFigFontNFSS{12}{14.4}{\rmdefault}{\mddefault}{\updefault}$\vertex=1\,,\lp=1$}}}}
\put(1420,78){\makebox(0,0)[lb]{\smash{{\SetFigFontNFSS{12}{14.4}{\rmdefault}{\mddefault}{\updefault}$\frac{1}{2}$}}}}
\end{picture}
}

%% file: figures/n3k0.eepic
\setlength{\unitlength}{0.00083333in}
\begingroup\makeatletter\ifx\SetFigFontNFSS\undefined%
\gdef\SetFigFontNFSS#1#2#3#4#5{%
  \reset@font\fontsize{#1}{#2pt}%
  \fontfamily{#3}\fontseries{#4}\fontshape{#5}%
  \selectfont}%
\fi\endgroup%
{\renewcommand{\dashlinestretch}{30}
\begin{picture}(2476,281)(0,-10)
\texture{44555555 55aaaaaa aa555555 55aaaaaa aa555555 55aaaaaa aa555555 55aaaaaa 
	aa555555 55aaaaaa aa555555 55aaaaaa aa555555 55aaaaaa aa555555 55aaaaaa 
	aa555555 55aaaaaa aa555555 55aaaaaa aa555555 55aaaaaa aa555555 55aaaaaa 
	aa555555 55aaaaaa aa555555 55aaaaaa aa555555 55aaaaaa aa555555 55aaaaaa }
\path(1790,158)(2090,158)
\path(1790,158)(2090,158)
\path(2124,161)(2424,161)
\path(2124,161)(2424,161)
\put(1792,158){\shade\ellipse{150}{150}}
\put(1792,158){\ellipse{150}{150}}
\put(2091,161){\shade\ellipse{150}{150}}
\put(2091,161){\ellipse{150}{150}}
\put(2393,166){\shade\ellipse{150}{150}}
\put(2393,166){\ellipse{150}{150}}
\put(1425,78){\makebox(0,0)[lb]{\smash{{\SetFigFontNFSS{12}{14.4}{\rmdefault}{\mddefault}{\updefault}$\frac{1}{2}$}}}}
\put(15,83){\makebox(0,0)[lb]{\smash{{\SetFigFontNFSS{12}{14.4}{\rmdefault}{\mddefault}{\updefault}$\vertex=3\,,\lp=0$}}}}
\end{picture}
}

%% file: figures/n2k1.eepic
\setlength{\unitlength}{0.00083333in}
\begingroup\makeatletter\ifx\SetFigFontNFSS\undefined%
\gdef\SetFigFontNFSS#1#2#3#4#5{%
  \reset@font\fontsize{#1}{#2pt}%
  \fontfamily{#3}\fontseries{#4}\fontshape{#5}%
  \selectfont}%
\fi\endgroup%
{\renewcommand{\dashlinestretch}{30}
\begin{picture}(3435,291)(0,-10)
\put(1805,136){\ellipse{240}{114}}
\path(1902,136)(2226,136)(2214,136)
\put(3203,159){\ellipse{300}{150}}
\texture{44555555 55aaaaaa aa555555 55aaaaaa aa555555 55aaaaaa aa555555 55aaaaaa 
	aa555555 55aaaaaa aa555555 55aaaaaa aa555555 55aaaaaa aa555555 55aaaaaa 
	aa555555 55aaaaaa aa555555 55aaaaaa aa555555 55aaaaaa aa555555 55aaaaaa 
	aa555555 55aaaaaa aa555555 55aaaaaa aa555555 55aaaaaa aa555555 55aaaaaa }
\put(2226,139){\shade\ellipse{150}{150}}
\put(2226,139){\ellipse{150}{150}}
\put(1916,135){\shade\ellipse{150}{150}}
\put(1916,135){\ellipse{150}{150}}
\put(3352,157){\shade\ellipse{150}{150}}
\put(3352,157){\ellipse{150}{150}}
\put(3047,156){\shade\ellipse{150}{150}}
\put(3047,156){\ellipse{150}{150}}
\put(15,93){\makebox(0,0)[lb]{\smash{{\SetFigFontNFSS{12}{14.4}{\rmdefault}{\mddefault}{\updefault}$\vertex=2\,,\lp=1$}}}}
\put(1425,78){\makebox(0,0)[lb]{\smash{{\SetFigFontNFSS{12}{14.4}{\rmdefault}{\mddefault}{\updefault}$\frac{1}{2}$}}}}
\put(2400,88){\makebox(0,0)[lb]{\smash{{\SetFigFontNFSS{12}{14.4}{\rmdefault}{\mddefault}{\updefault}$+\,\,\frac{1}{2^ 2}$}}}}
\end{picture}
}

%% file: figures/n1k2.eepic
\setlength{\unitlength}{0.00083333in}
\begingroup\makeatletter\ifx\SetFigFontNFSS\undefined%
\gdef\SetFigFontNFSS#1#2#3#4#5{%
  \reset@font\fontsize{#1}{#2pt}%
  \fontfamily{#3}\fontseries{#4}\fontshape{#5}%
  \selectfont}%
\fi\endgroup%
{\renewcommand{\dashlinestretch}{30}
\begin{picture}(1936,528)(0,-10)
\put(1853,384){\ellipse{120}{246}}
\put(1853,129){\ellipse{120}{246}}
\texture{44555555 55aaaaaa aa555555 55aaaaaa aa555555 55aaaaaa aa555555 55aaaaaa 
	aa555555 55aaaaaa aa555555 55aaaaaa aa555555 55aaaaaa aa555555 55aaaaaa 
	aa555555 55aaaaaa aa555555 55aaaaaa aa555555 55aaaaaa aa555555 55aaaaaa 
	aa555555 55aaaaaa aa555555 55aaaaaa aa555555 55aaaaaa aa555555 55aaaaaa }
\put(1853,257){\shade\ellipse{150}{150}}
\put(1853,257){\ellipse{150}{150}}
\put(15,204){\makebox(0,0)[lb]{\smash{{\SetFigFontNFSS{12}{14.4}{\rmdefault}{\mddefault}{\updefault}$\vertex=1\,,\lp=2$}}}}
\put(1440,209){\makebox(0,0)[lb]{\smash{{\SetFigFontNFSS{12}{14.4}{\rmdefault}{\mddefault}{\updefault}$\frac{1}{2^3}$}}}}
\end{picture}
}

%% file: figures/n4k0.eepic
\setlength{\unitlength}{0.00083333in}
\begingroup\makeatletter\ifx\SetFigFontNFSS\undefined%
\gdef\SetFigFontNFSS#1#2#3#4#5{%
  \reset@font\fontsize{#1}{#2pt}%
  \fontfamily{#3}\fontseries{#4}\fontshape{#5}%
  \selectfont}%
\fi\endgroup%
{\renewcommand{\dashlinestretch}{30}
\begin{picture}(4018,701)(0,-10)
\texture{44555555 55aaaaaa aa555555 55aaaaaa aa555555 55aaaaaa aa555555 55aaaaaa 
	aa555555 55aaaaaa aa555555 55aaaaaa aa555555 55aaaaaa aa555555 55aaaaaa 
	aa555555 55aaaaaa aa555555 55aaaaaa aa555555 55aaaaaa aa555555 55aaaaaa 
	aa555555 55aaaaaa aa555555 55aaaaaa aa555555 55aaaaaa aa555555 55aaaaaa }
\path(1831,331)(2131,331)
\path(1831,331)(2131,331)
\path(2131,331)(2431,331)
\path(2131,331)(2431,331)
\path(2431,331)(2731,331)
\path(2431,331)(2731,331)
\path(3935,343)(3635,343)
\path(3935,343)(3635,343)
\path(3485,603)(3645,333)
\path(3495,93)(3655,363)
\put(2101,337){\shade\ellipse{150}{150}}
\put(2101,337){\ellipse{150}{150}}
\put(2401,351){\shade\ellipse{150}{150}}
\put(2401,351){\ellipse{150}{150}}
\put(2694,344){\shade\ellipse{150}{150}}
\put(2694,344){\ellipse{150}{150}}
\put(1798,345){\shade\ellipse{150}{150}}
\put(1798,345){\ellipse{150}{150}}
\put(3635,343){\shade\ellipse{150}{150}}
\put(3635,343){\ellipse{150}{150}}
\put(3935,343){\shade\ellipse{150}{150}}
\put(3935,343){\ellipse{150}{150}}
\put(3490,603){\shade\ellipse{150}{150}}
\put(3490,603){\ellipse{150}{150}}
\put(3490,83){\shade\ellipse{150}{150}}
\put(3490,83){\ellipse{150}{150}}
\put(1455,278){\makebox(0,0)[lb]{\smash{{\SetFigFontNFSS{12}{14.4}{\rmdefault}{\mddefault}{\updefault}$\frac{1}{2}$}}}}
\put(15,275){\makebox(0,0)[lb]{\smash{{\SetFigFontNFSS{12}{14.4}{\rmdefault}{\mddefault}{\updefault}$\vertex=4\,,\lp=0$}}}}
\put(2893,280){\makebox(0,0)[lb]{\smash{{\SetFigFontNFSS{12}{14.4}{\rmdefault}{\mddefault}{\updefault}$+\,\,\frac{1}{3!}$}}}}
\end{picture}
}

%% file: figures/n3k1.eepic
\setlength{\unitlength}{0.00083333in}
\begingroup\makeatletter\ifx\SetFigFontNFSS\undefined%
\gdef\SetFigFontNFSS#1#2#3#4#5{%
  \reset@font\fontsize{#1}{#2pt}%
  \fontfamily{#3}\fontseries{#4}\fontshape{#5}%
  \selectfont}%
\fi\endgroup%
{\renewcommand{\dashlinestretch}{30}
\begin{picture}(6547,507)(0,-10)
\put(3298,290){\ellipse{120}{246}}
\put(4398,163){\ellipse{240}{114}}
\path(2940,166)(3300,166)(3306,166)
\path(3240,166)(3600,166)(3606,166)
\path(5200,170)(4840,170)(4834,170)
\path(4900,170)(4540,170)(4534,170)
\put(6031,178){\ellipse{300}{150}}
\path(6140,175)(6500,175)(6506,175)
\path(1982,403)(1839,155)(2125,155)(1982,403)
\texture{44555555 55aaaaaa aa555555 55aaaaaa aa555555 55aaaaaa aa555555 55aaaaaa 
	aa555555 55aaaaaa aa555555 55aaaaaa aa555555 55aaaaaa aa555555 55aaaaaa 
	aa555555 55aaaaaa aa555555 55aaaaaa aa555555 55aaaaaa aa555555 55aaaaaa 
	aa555555 55aaaaaa aa555555 55aaaaaa aa555555 55aaaaaa aa555555 55aaaaaa }
\put(1979,411){\shade\ellipse{150}{150}}
\put(1979,411){\ellipse{150}{150}}
\put(1832,155){\shade\ellipse{150}{150}}
\put(1832,155){\ellipse{150}{150}}
\put(2135,162){\shade\ellipse{150}{150}}
\put(2135,162){\ellipse{150}{150}}
\put(3300,171){\shade\ellipse{150}{150}}
\put(3300,171){\ellipse{150}{150}}
\put(3588,177){\shade\ellipse{150}{150}}
\put(3588,177){\ellipse{150}{150}}
\put(3008,172){\shade\ellipse{150}{150}}
\put(3008,172){\ellipse{150}{150}}
\put(4840,165){\shade\ellipse{150}{150}}
\put(4840,165){\ellipse{150}{150}}
\put(5132,164){\shade\ellipse{150}{150}}
\put(5132,164){\ellipse{150}{150}}
\put(4534,168){\shade\ellipse{150}{150}}
\put(4534,168){\ellipse{150}{150}}
\put(5882,178){\shade\ellipse{150}{150}}
\put(5882,178){\ellipse{150}{150}}
\put(6176,178){\shade\ellipse{150}{150}}
\put(6176,178){\ellipse{150}{150}}
\put(6464,184){\shade\ellipse{150}{150}}
\put(6464,184){\ellipse{150}{150}}
\put(1450,78){\makebox(0,0)[lb]{\smash{{\SetFigFontNFSS{12}{14.4}{\rmdefault}{\mddefault}{\updefault}$\frac{1}{3!}$}}}}
\put(2370,83){\makebox(0,0)[lb]{\smash{{\SetFigFontNFSS{12}{14.4}{\rmdefault}{\mddefault}{\updefault}$+\,\,\frac{1}{2^2}$}}}}
\put(3820,103){\makebox(0,0)[lb]{\smash{{\SetFigFontNFSS{12}{14.4}{\rmdefault}{\mddefault}{\updefault}$+\,\,\frac{1}{2}$}}}}
\put(5320,118){\makebox(0,0)[lb]{\smash{{\SetFigFontNFSS{12}{14.4}{\rmdefault}{\mddefault}{\updefault}$+\,\,\frac{1}{2}$}}}}
\put(15,78){\makebox(0,0)[lb]{\smash{{\SetFigFontNFSS{12}{14.4}{\rmdefault}{\mddefault}{\updefault}$\vertex=3\,,\lp=1$}}}}
\end{picture}
}

%% file: figures/n2k2.eepic
\setlength{\unitlength}{0.00083333in}
\begingroup\makeatletter\ifx\SetFigFontNFSS\undefined%
\gdef\SetFigFontNFSS#1#2#3#4#5{%
  \reset@font\fontsize{#1}{#2pt}%
  \fontfamily{#3}\fontseries{#4}\fontshape{#5}%
  \selectfont}%
\fi\endgroup%
{\renewcommand{\dashlinestretch}{30}
\begin{picture}(6499,507)(0,-10)
\put(1896,363){\ellipse{120}{246}}
\put(1896,129){\ellipse{120}{246}}
\put(3583,260){\ellipse{240}{114}}
\put(3046,256){\ellipse{240}{114}}
\put(4575,243){\ellipse{240}{114}}
\path(1878,244)(2202,244)(2190,244)
\path(3486,260)(3162,260)(3174,260)
\put(4845,245){\ellipse{300}{150}}
\put(6270,270){\ellipse{300}{150}}
\path(6092,267)(6416,267)(6404,267)
\texture{44555555 55aaaaaa aa555555 55aaaaaa aa555555 55aaaaaa aa555555 55aaaaaa 
	aa555555 55aaaaaa aa555555 55aaaaaa aa555555 55aaaaaa aa555555 55aaaaaa 
	aa555555 55aaaaaa aa555555 55aaaaaa aa555555 55aaaaaa aa555555 55aaaaaa 
	aa555555 55aaaaaa aa555555 55aaaaaa aa555555 55aaaaaa aa555555 55aaaaaa }
\put(2202,247){\shade\ellipse{150}{150}}
\put(2202,247){\ellipse{150}{150}}
\put(1892,243){\shade\ellipse{150}{150}}
\put(1892,243){\ellipse{150}{150}}
\put(3162,257){\shade\ellipse{150}{150}}
\put(3162,257){\ellipse{150}{150}}
\put(3472,261){\shade\ellipse{150}{150}}
\put(3472,261){\ellipse{150}{150}}
\put(4993,245){\shade\ellipse{150}{150}}
\put(4993,245){\ellipse{150}{150}}
\put(4684,241){\shade\ellipse{150}{150}}
\put(4684,241){\ellipse{150}{150}}
\put(6416,270){\shade\ellipse{150}{150}}
\put(6416,270){\ellipse{150}{150}}
\put(6122,270){\shade\ellipse{150}{150}}
\put(6122,270){\ellipse{150}{150}}
\put(15,176){\makebox(0,0)[lb]{\smash{{\SetFigFontNFSS{12}{14.4}{\rmdefault}{\mddefault}{\updefault}$\vertex=2\,,\lp=2$}}}}
\put(1445,181){\makebox(0,0)[lb]{\smash{{\SetFigFontNFSS{12}{14.4}{\rmdefault}{\mddefault}{\updefault}$\frac{1}{2^3}$}}}}
\put(2390,176){\makebox(0,0)[lb]{\smash{{\SetFigFontNFSS{12}{14.4}{\rmdefault}{\mddefault}{\updefault}$+\,\,\frac{1}{2^ 3}$}}}}
\put(5365,201){\makebox(0,0)[lb]{\smash{{\SetFigFontNFSS{12}{14.4}{\rmdefault}{\mddefault}{\updefault}$+\,\,\frac{1}{2\cdot3!}$}}}}
\put(3865,201){\makebox(0,0)[lb]{\smash{{\SetFigFontNFSS{12}{14.4}{\rmdefault}{\mddefault}{\updefault}$+\,\,\frac{1}{2^ 2\,\,}$}}}}
\end{picture}
}

%% file: figures/n1k3.eepic
\setlength{\unitlength}{0.00083333in}
\begingroup\makeatletter\ifx\SetFigFontNFSS\undefined%
\gdef\SetFigFontNFSS#1#2#3#4#5{%
  \reset@font\fontsize{#1}{#2pt}%
  \fontfamily{#3}\fontseries{#4}\fontshape{#5}%
  \selectfont}%
\fi\endgroup%
{\renewcommand{\dashlinestretch}{30}
\begin{picture}(2416,386)(0,-10)
\put(2163,241){\ellipse{120}{246}}
\put(2048,112){\ellipse{240}{114}}
\put(2288,112){\ellipse{240}{114}}
\texture{44555555 55aaaaaa aa555555 55aaaaaa aa555555 55aaaaaa aa555555 55aaaaaa 
	aa555555 55aaaaaa aa555555 55aaaaaa aa555555 55aaaaaa aa555555 55aaaaaa 
	aa555555 55aaaaaa aa555555 55aaaaaa aa555555 55aaaaaa aa555555 55aaaaaa 
	aa555555 55aaaaaa aa555555 55aaaaaa aa555555 55aaaaaa aa555555 55aaaaaa }
\put(2168,109){\shade\ellipse{150}{150}}
\put(2168,109){\ellipse{150}{150}}
\put(1435,78){\makebox(0,0)[lb]{\smash{{\SetFigFontNFSS{12}{14.4}{\rmdefault}{\mddefault}{\updefault}$\frac{1}{2^ 3\cdot3!}$}}}}
\put(15,83){\makebox(0,0)[lb]{\smash{{\SetFigFontNFSS{12}{14.4}{\rmdefault}{\mddefault}{\updefault}$\vertex=1\,,\lp=3$}}}}
\end{picture}
}

%% file: appbiconn.tex
\section{}\label{app:biconn}
This appendix  shows  $\Dp^{\vertex,\lp,0}$ up to order $\vertex+\lp\le5$ and computed via formula  (\ref{eq:recbiconn}).  All graphs in the same equivalence class are identified as the same.   The coefficients in front  of graphs are the
inverses of the orders of their
groups of automorphisms.
\vspace{1cm}\\
\input{figures/n1k0.eepic}
\vspace{1cm}\\
\input{figures/n1k1.eepic}
\vspace{1cm}\\
\input{figures/noloopsn2k1.eepic}
\vspace{1cm}\\
\input{figures/n1k2.eepic}
\vspace{1cm}\\
\input{figures/bin3k1.eepic}
\vspace{1cm}\\
\input{figures/bin2k2.eepic}
\vspace{1cm}\\
\input{figures/n1k3.eepic}
\vspace{1cm}\\
\input{figures/biconnn4k1.eepic}
\vspace{1cm}\\
\input{figures/bin3k2.eepic}
\vspace{1cm}\\
\input{figures/bin2k3.eepic}
\vspace{0,6cm}\\
\input{figures/n1k4.eepic}

%% file: figures/noloopsn2k1.eepic
\setlength{\unitlength}{0.00083333in}
\begingroup\makeatletter\ifx\SetFigFontNFSS\undefined%
\gdef\SetFigFontNFSS#1#2#3#4#5{%
  \reset@font\fontsize{#1}{#2pt}%
  \fontfamily{#3}\fontseries{#4}\fontshape{#5}%
  \selectfont}%
\fi\endgroup%
{\renewcommand{\dashlinestretch}{30}
\begin{picture}(2231,281)(0,-10)
\put(1999,159){\ellipse{300}{150}}
\texture{44555555 55aaaaaa aa555555 55aaaaaa aa555555 55aaaaaa aa555555 55aaaaaa 
	aa555555 55aaaaaa aa555555 55aaaaaa aa555555 55aaaaaa aa555555 55aaaaaa 
	aa555555 55aaaaaa aa555555 55aaaaaa aa555555 55aaaaaa aa555555 55aaaaaa 
	aa555555 55aaaaaa aa555555 55aaaaaa aa555555 55aaaaaa aa555555 55aaaaaa }
\put(2148,157){\shade\ellipse{150}{150}}
\put(2148,157){\ellipse{150}{150}}
\put(1843,156){\shade\ellipse{150}{150}}
\put(1843,156){\ellipse{150}{150}}
\put(15,83){\makebox(0,0)[lb]{\smash{{\SetFigFontNFSS{12}{14.4}{\rmdefault}{\mddefault}{\updefault}$\vertex=2\,,\lp=1$}}}}
\put(1425,78){\makebox(0,0)[lb]{\smash{{\SetFigFontNFSS{12}{14.4}{\rmdefault}{\mddefault}{\updefault}$\frac{1}{2^ 2}$}}}}
\end{picture}
}

%% file: figures/bin3k1.eepic
\setlength{\unitlength}{0.00083333in}
\begingroup\makeatletter\ifx\SetFigFontNFSS\undefined%
\gdef\SetFigFontNFSS#1#2#3#4#5{%
  \reset@font\fontsize{#1}{#2pt}%
  \fontfamily{#3}\fontseries{#4}\fontshape{#5}%
  \selectfont}%
\fi\endgroup%
{\renewcommand{\dashlinestretch}{30}
\begin{picture}(2220,487)(0,-10)
\path(1984,383)(1841,135)(2127,135)(1984,383)
\texture{44555555 55aaaaaa aa555555 55aaaaaa aa555555 55aaaaaa aa555555 55aaaaaa 
	aa555555 55aaaaaa aa555555 55aaaaaa aa555555 55aaaaaa aa555555 55aaaaaa 
	aa555555 55aaaaaa aa555555 55aaaaaa aa555555 55aaaaaa aa555555 55aaaaaa 
	aa555555 55aaaaaa aa555555 55aaaaaa aa555555 55aaaaaa aa555555 55aaaaaa }
\put(1981,391){\shade\ellipse{150}{150}}
\put(1981,391){\ellipse{150}{150}}
\put(1834,135){\shade\ellipse{150}{150}}
\put(1834,135){\ellipse{150}{150}}
\put(2137,142){\shade\ellipse{150}{150}}
\put(2137,142){\ellipse{150}{150}}
\put(15,83){\makebox(0,0)[lb]{\smash{{\SetFigFontNFSS{12}{14.4}{\rmdefault}{\mddefault}{\updefault}$\vertex=3\,,\lp=1$}}}}
\put(1435,78){\makebox(0,0)[lb]{\smash{{\SetFigFontNFSS{12}{14.4}{\rmdefault}{\mddefault}{\updefault}$\frac{1}{3!}$}}}}
\end{picture}
}

%% file: figures/bin2k2.eepic
\setlength{\unitlength}{0.00083333in}
\begingroup\makeatletter\ifx\SetFigFontNFSS\undefined%
\gdef\SetFigFontNFSS#1#2#3#4#5{%
  \reset@font\fontsize{#1}{#2pt}%
  \fontfamily{#3}\fontseries{#4}\fontshape{#5}%
  \selectfont}%
\fi\endgroup%
{\renewcommand{\dashlinestretch}{30}
\begin{picture}(3490,281)(0,-10)
\put(1860,115){\ellipse{240}{114}}
\put(3261,132){\ellipse{300}{150}}
\path(3083,129)(3407,129)(3395,129)
\put(2130,117){\ellipse{300}{150}}
\texture{44555555 55aaaaaa aa555555 55aaaaaa aa555555 55aaaaaa aa555555 55aaaaaa 
	aa555555 55aaaaaa aa555555 55aaaaaa aa555555 55aaaaaa aa555555 55aaaaaa 
	aa555555 55aaaaaa aa555555 55aaaaaa aa555555 55aaaaaa aa555555 55aaaaaa 
	aa555555 55aaaaaa aa555555 55aaaaaa aa555555 55aaaaaa aa555555 55aaaaaa }
\put(3407,132){\shade\ellipse{150}{150}}
\put(3407,132){\ellipse{150}{150}}
\put(3113,132){\shade\ellipse{150}{150}}
\put(3113,132){\ellipse{150}{150}}
\put(2278,117){\shade\ellipse{150}{150}}
\put(2278,117){\ellipse{150}{150}}
\put(1969,113){\shade\ellipse{150}{150}}
\put(1969,113){\ellipse{150}{150}}
\put(15,83){\makebox(0,0)[lb]{\smash{{\SetFigFontNFSS{12}{14.4}{\rmdefault}{\mddefault}{\updefault}$\vertex=2\,,\lp=2$}}}}
\put(2510,83){\makebox(0,0)[lb]{\smash{{\SetFigFontNFSS{12}{14.4}{\rmdefault}{\mddefault}{\updefault}$+\,\,\frac{1}{2\cdot3!}$}}}}
\put(1420,78){\makebox(0,0)[lb]{\smash{{\SetFigFontNFSS{12}{14.4}{\rmdefault}{\mddefault}{\updefault}$\frac{1}{2^ 2}$}}}}
\end{picture}
}

%% file: figures/biconnn4k1.eepic
\setlength{\unitlength}{0.00083333in}
\begingroup\makeatletter\ifx\SetFigFontNFSS\undefined%
\gdef\SetFigFontNFSS#1#2#3#4#5{%
  \reset@font\fontsize{#1}{#2pt}%
  \fontfamily{#3}\fontseries{#4}\fontshape{#5}%
  \selectfont}%
\fi\endgroup%
{\renewcommand{\dashlinestretch}{30}
\begin{picture}(2231,482)(0,-10)
\path(2149,423)(2149,25)
\path(1843,429)(1843,31)
\texture{44555555 55aaaaaa aa555555 55aaaaaa aa555555 55aaaaaa aa555555 55aaaaaa 
	aa555555 55aaaaaa aa555555 55aaaaaa aa555555 55aaaaaa aa555555 55aaaaaa 
	aa555555 55aaaaaa aa555555 55aaaaaa aa555555 55aaaaaa aa555555 55aaaaaa 
	aa555555 55aaaaaa aa555555 55aaaaaa aa555555 55aaaaaa aa555555 55aaaaaa }
\path(1861,84)(2161,84)
\path(1861,84)(2161,84)
\path(1861,385)(2161,385)
\path(1861,385)(2161,385)
\put(2148,380){\shade\ellipse{150}{150}}
\put(2148,380){\ellipse{150}{150}}
\put(2148,85){\shade\ellipse{150}{150}}
\put(2148,85){\ellipse{150}{150}}
\put(1841,384){\shade\ellipse{150}{150}}
\put(1841,384){\ellipse{150}{150}}
\put(1850,83){\shade\ellipse{150}{150}}
\put(1850,83){\ellipse{150}{150}}
\put(15,150){\makebox(0,0)[lb]{\smash{{\SetFigFontNFSS{12}{14.4}{\rmdefault}{\mddefault}{\updefault}$\vertex=4\,,\lp=1$}}}}
\put(1425,150){\makebox(0,0)[lb]{\smash{{\SetFigFontNFSS{12}{14.4}{\rmdefault}{\mddefault}{\updefault}$\frac{1}{8}$}}}}
\end{picture}
}

%% file: figures/bin3k2.eepic
\setlength{\unitlength}{0.00083333in}
\begingroup\makeatletter\ifx\SetFigFontNFSS\undefined%
\gdef\SetFigFontNFSS#1#2#3#4#5{%
  \reset@font\fontsize{#1}{#2pt}%
  \fontfamily{#3}\fontseries{#4}\fontshape{#5}%
  \selectfont}%
\fi\endgroup%
{\renewcommand{\dashlinestretch}{30}
\begin{picture}(4839,607)(0,-10)
\path(3006,93)(3153,466)
\path(3313,99)(3159,466)
\drawline(3159,466)(3159,466)
\put(1959,463){\ellipse{120}{246}}
\put(3166,93){\ellipse{300}{150}}
\put(4306,186){\ellipse{300}{150}}
\put(4608,180){\ellipse{300}{150}}
\path(1960,331)(1817,83)(2103,83)(1960,331)
\texture{44555555 55aaaaaa aa555555 55aaaaaa aa555555 55aaaaaa aa555555 55aaaaaa 
	aa555555 55aaaaaa aa555555 55aaaaaa aa555555 55aaaaaa aa555555 55aaaaaa 
	aa555555 55aaaaaa aa555555 55aaaaaa aa555555 55aaaaaa aa555555 55aaaaaa 
	aa555555 55aaaaaa aa555555 55aaaaaa aa555555 55aaaaaa aa555555 55aaaaaa }
\put(1810,83){\shade\ellipse{150}{150}}
\put(1810,83){\ellipse{150}{150}}
\put(2113,90){\shade\ellipse{150}{150}}
\put(2113,90){\ellipse{150}{150}}
\put(1957,334){\shade\ellipse{150}{150}}
\put(1957,334){\ellipse{150}{150}}
\put(3010,90){\shade\ellipse{150}{150}}
\put(3010,90){\ellipse{150}{150}}
\put(3310,96){\shade\ellipse{150}{150}}
\put(3310,96){\ellipse{150}{150}}
\put(3157,462){\shade\ellipse{150}{150}}
\put(3157,462){\ellipse{150}{150}}
\put(4157,186){\shade\ellipse{150}{150}}
\put(4157,186){\ellipse{150}{150}}
\put(4451,186){\shade\ellipse{150}{150}}
\put(4451,186){\ellipse{150}{150}}
\put(4756,187){\shade\ellipse{150}{150}}
\put(4756,187){\ellipse{150}{150}}
\put(15,96){\makebox(0,0)[lb]{\smash{{\SetFigFontNFSS{12}{14.4}{\rmdefault}{\mddefault}{\updefault}$\vertex=3\,,\lp=2$}}}}
\put(1426,95){\makebox(0,0)[lb]{\smash{{\SetFigFontNFSS{12}{14.4}{\rmdefault}{\mddefault}{\updefault}$\frac{1}{2^2}$}}}}
\put(2366,95){\makebox(0,0)[lb]{\smash{{\SetFigFontNFSS{12}{14.4}{\rmdefault}{\mddefault}{\updefault}$+\,\,\frac{1}{2^2}$}}}}
\put(3541,100){\makebox(0,0)[lb]{\smash{{\SetFigFontNFSS{12}{14.4}{\rmdefault}{\mddefault}{\updefault}$+\,\,\frac{1}{2^3}$}}}}
\end{picture}
}

%% file: figures/bin2k3.eepic
\setlength{\unitlength}{0.00083333in}
\begingroup\makeatletter\ifx\SetFigFontNFSS\undefined%
\gdef\SetFigFontNFSS#1#2#3#4#5{%
  \reset@font\fontsize{#1}{#2pt}%
  \fontfamily{#3}\fontseries{#4}\fontshape{#5}%
  \selectfont}%
\fi\endgroup%
{\renewcommand{\dashlinestretch}{30}
\begin{picture}(6578,524)(0,-10)
\put(6352,269){\ellipse{310}{310}}
\put(2488,233){\ellipse{240}{114}}
\put(1956,223){\ellipse{240}{114}}
\put(3607,260){\ellipse{240}{114}}
\put(4941,130){\ellipse{120}{246}}
\put(4943,380){\ellipse{120}{246}}
\put(2217,229){\ellipse{300}{150}}
\put(3872,259){\ellipse{300}{150}}
\path(3694,256)(4018,256)(4006,256)
\put(5091,256){\ellipse{300}{150}}
\put(6347,267){\ellipse{300}{150}}
\texture{44555555 55aaaaaa aa555555 55aaaaaa aa555555 55aaaaaa aa555555 55aaaaaa 
	aa555555 55aaaaaa aa555555 55aaaaaa aa555555 55aaaaaa aa555555 55aaaaaa 
	aa555555 55aaaaaa aa555555 55aaaaaa aa555555 55aaaaaa aa555555 55aaaaaa 
	aa555555 55aaaaaa aa555555 55aaaaaa aa555555 55aaaaaa aa555555 55aaaaaa }
\put(2369,228){\shade\ellipse{150}{150}}
\put(2369,228){\ellipse{150}{150}}
\put(2075,228){\shade\ellipse{150}{150}}
\put(2075,228){\ellipse{150}{150}}
\put(4018,259){\shade\ellipse{150}{150}}
\put(4018,259){\ellipse{150}{150}}
\put(3724,259){\shade\ellipse{150}{150}}
\put(3724,259){\ellipse{150}{150}}
\put(5237,253){\shade\ellipse{150}{150}}
\put(5237,253){\ellipse{150}{150}}
\put(4943,253){\shade\ellipse{150}{150}}
\put(4943,253){\ellipse{150}{150}}
\put(6203,270){\shade\ellipse{150}{150}}
\put(6203,270){\ellipse{150}{150}}
\put(6495,273){\shade\ellipse{150}{150}}
\put(6495,273){\ellipse{150}{150}}
\put(15,147){\makebox(0,0)[lb]{\smash{{\SetFigFontNFSS{12}{14.4}{\rmdefault}{\mddefault}{\updefault}$\vertex=2\,,\lp=3$}}}}
\put(2841,165){\makebox(0,0)[lb]{\smash{{\SetFigFontNFSS{12}{14.4}{\rmdefault}{\mddefault}{\updefault}$+\,\,\,\frac{1}{2\cdot3!}$}}}}
\put(4261,165){\makebox(0,0)[lb]{\smash{{\SetFigFontNFSS{12}{14.4}{\rmdefault}{\mddefault}{\updefault}$+\,\,\,\frac{1}{2^ 4}$}}}}
\put(5491,170){\makebox(0,0)[lb]{\smash{{\SetFigFontNFSS{12}{14.4}{\rmdefault}{\mddefault}{\updefault}$+\,\,\,\frac{1}{2\cdot4!}$}}}}
\put(1496,165){\makebox(0,0)[lb]{\smash{{\SetFigFontNFSS{12}{14.4}{\rmdefault}{\mddefault}{\updefault}$\frac{1}{2^ 4}$}}}}
\end{picture}
}

%% file: figures/n1k4.eepic
\setlength{\unitlength}{0.00083333in}
\begingroup\makeatletter\ifx\SetFigFontNFSS\undefined%
\gdef\SetFigFontNFSS#1#2#3#4#5{%
  \reset@font\fontsize{#1}{#2pt}%
  \fontfamily{#3}\fontseries{#4}\fontshape{#5}%
  \selectfont}%
\fi\endgroup%
{\renewcommand{\dashlinestretch}{30}
\begin{picture}(2491,512)(0,-10)
\put(2123,230){\ellipse{240}{114}}
\put(2363,230){\ellipse{240}{114}}
\put(2239,130){\ellipse{120}{246}}
\put(2248,367){\ellipse{120}{246}}
\texture{44555555 55aaaaaa aa555555 55aaaaaa aa555555 55aaaaaa aa555555 55aaaaaa 
	aa555555 55aaaaaa aa555555 55aaaaaa aa555555 55aaaaaa aa555555 55aaaaaa 
	aa555555 55aaaaaa aa555555 55aaaaaa aa555555 55aaaaaa aa555555 55aaaaaa 
	aa555555 55aaaaaa aa555555 55aaaaaa aa555555 55aaaaaa aa555555 55aaaaaa }
\put(2243,227){\shade\ellipse{150}{150}}
\put(2243,227){\ellipse{150}{150}}
\put(15,166){\makebox(0,0)[lb]{\smash{{\SetFigFontNFSS{12}{14.4}{\rmdefault}{\mddefault}{\updefault}$\vertex=1\,,\lp=4$}}}}
\put(1485,166){\makebox(0,0)[lb]{\smash{{\SetFigFontNFSS{12}{14.4}{\rmdefault}{\mddefault}{\updefault}$\frac{1}{2^ 4\cdot4!}$}}}}
\end{picture}
}

%% file: appsimple.tex
\section{}\label{app:simple}
This appendix shows  $\Dnew^{\vertex,\lp,0}$ up to order $\vertex+\lp\le6$  and computed via formula  (\ref{eq:recsimple}).  All graphs in the same equivalence class are identified as the same.   The coefficients in front  of graphs are the
inverses of the orders of their
groups of automorphisms.
\vspace{1cm}\\
\input{figures/n1k0.eepic}
\vspace{0.95cm}\\
\input{figures/n2k0.eepic}
\vspace{0.95cm}\\
\input{figures/n3k0.eepic}
\vspace{0.9cm}\\
\input{figures/n4k0.eepic}
\vspace{0.5cm}\\
\input{figures/bin3k1.eepic}
\vspace{0.9cm}\\
\input{figures/n5k0.eepic}
\vspace{1cm}\\
\input{figures/simplen4k1.eepic}
\vspace{0.9cm}\\
\input{figures/n6k0.eepic}
\vspace{1cm}\\
\input{figures/simplen5k1.eepic}
\vspace{1cm}\\
\input{figures/simplen4k2.eepic}

%% file: figures/n5k0.eepic
\setlength{\unitlength}{0.00083333in}
\begingroup\makeatletter\ifx\SetFigFontNFSS\undefined%
\gdef\SetFigFontNFSS#1#2#3#4#5{%
  \reset@font\fontsize{#1}{#2pt}%
  \fontfamily{#3}\fontseries{#4}\fontshape{#5}%
  \selectfont}%
\fi\endgroup%
{\renewcommand{\dashlinestretch}{30}
\begin{picture}(5840,781)(0,-10)
\texture{44555555 55aaaaaa aa555555 55aaaaaa aa555555 55aaaaaa aa555555 55aaaaaa 
	aa555555 55aaaaaa aa555555 55aaaaaa aa555555 55aaaaaa aa555555 55aaaaaa 
	aa555555 55aaaaaa aa555555 55aaaaaa aa555555 55aaaaaa aa555555 55aaaaaa 
	aa555555 55aaaaaa aa555555 55aaaaaa aa555555 55aaaaaa aa555555 55aaaaaa }
\path(1832,371)(2132,371)
\path(1832,371)(2132,371)
\path(2132,371)(2432,371)
\path(2132,371)(2432,371)
\path(2432,371)(2732,371)
\path(2432,371)(2732,371)
\path(2732,371)(3032,371)
\path(2732,371)(3032,371)
\path(4045,378)(3745,378)
\path(4045,378)(3745,378)
\path(3595,638)(3755,368)
\path(3605,128)(3765,398)
\path(4315,378)(4015,378)
\path(4315,378)(4015,378)
\path(5457,83)(5457,383)
\path(5457,83)(5457,383)
\path(5457,383)(5457,683)
\path(5457,383)(5457,683)
\path(5532,383)(5232,383)
\path(5532,383)(5232,383)
\path(5457,83)(5457,383)
\path(5457,83)(5457,383)
\path(5757,383)(5457,383)
\path(5757,383)(5457,383)
\put(1832,371){\shade\ellipse{150}{150}}
\put(1832,371){\ellipse{150}{150}}
\put(2132,371){\shade\ellipse{150}{150}}
\put(2132,371){\ellipse{150}{150}}
\put(2432,371){\shade\ellipse{150}{150}}
\put(2432,371){\ellipse{150}{150}}
\put(2732,371){\shade\ellipse{150}{150}}
\put(2732,371){\ellipse{150}{150}}
\put(3032,371){\shade\ellipse{150}{150}}
\put(3032,371){\ellipse{150}{150}}
\put(3745,378){\shade\ellipse{150}{150}}
\put(3745,378){\ellipse{150}{150}}
\put(4045,378){\shade\ellipse{150}{150}}
\put(4045,378){\ellipse{150}{150}}
\put(4345,378){\shade\ellipse{150}{150}}
\put(4345,378){\ellipse{150}{150}}
\put(3600,638){\shade\ellipse{150}{150}}
\put(3600,638){\ellipse{150}{150}}
\put(3600,118){\shade\ellipse{150}{150}}
\put(3600,118){\ellipse{150}{150}}
\put(5157,383){\shade\ellipse{150}{150}}
\put(5157,383){\ellipse{150}{150}}
\put(5457,83){\shade\ellipse{150}{150}}
\put(5457,83){\ellipse{150}{150}}
\put(5457,383){\shade\ellipse{150}{150}}
\put(5457,383){\ellipse{150}{150}}
\put(5457,683){\shade\ellipse{150}{150}}
\put(5457,683){\ellipse{150}{150}}
\put(5757,383){\shade\ellipse{150}{150}}
\put(5757,383){\ellipse{150}{150}}
\put(1445,300){\makebox(0,0)[lb]{\smash{{\SetFigFontNFSS{12}{14.4}{\rmdefault}{\mddefault}{\updefault}$\frac{1}{2}$}}}}
\put(15,310){\makebox(0,0)[lb]{\smash{{\SetFigFontNFSS{12}{14.4}{\rmdefault}{\mddefault}{\updefault}$\vertex=5\,,\lp=0$}}}}
\put(3195,313){\makebox(0,0)[lb]{\smash{{\SetFigFontNFSS{12}{14.4}{\rmdefault}{\mddefault}{\updefault}$+\,\,\frac{1}{2}$}}}}
\put(4570,313){\makebox(0,0)[lb]{\smash{{\SetFigFontNFSS{12}{14.4}{\rmdefault}{\mddefault}{\updefault}$+\,\,\frac{1}{4!}$}}}}
\end{picture}
}

%% file: figures/simplen4k1.eepic
\setlength{\unitlength}{0.00083333in}
\begingroup\makeatletter\ifx\SetFigFontNFSS\undefined%
\gdef\SetFigFontNFSS#1#2#3#4#5{%
  \reset@font\fontsize{#1}{#2pt}%
  \fontfamily{#3}\fontseries{#4}\fontshape{#5}%
  \selectfont}%
\fi\endgroup%
{\renewcommand{\dashlinestretch}{30}
\begin{picture}(3705,496)(0,-10)
\path(2189,437)(2189,39)
\path(1883,443)(1883,45)
\texture{44555555 55aaaaaa aa555555 55aaaaaa aa555555 55aaaaaa aa555555 55aaaaaa 
	aa555555 55aaaaaa aa555555 55aaaaaa aa555555 55aaaaaa aa555555 55aaaaaa 
	aa555555 55aaaaaa aa555555 55aaaaaa aa555555 55aaaaaa aa555555 55aaaaaa 
	aa555555 55aaaaaa aa555555 55aaaaaa aa555555 55aaaaaa aa555555 55aaaaaa }
\path(1901,98)(2201,98)
\path(1901,98)(2201,98)
\path(1901,399)(2201,399)
\path(1901,399)(2201,399)
\path(3317,236)(3069,379)(3069,93)(3317,236)
\path(3390,240)(3690,240)
\path(3390,240)(3690,240)
\put(2188,394){\shade\ellipse{150}{150}}
\put(2188,394){\ellipse{150}{150}}
\put(2188,99){\shade\ellipse{150}{150}}
\put(2188,99){\ellipse{150}{150}}
\put(1881,398){\shade\ellipse{150}{150}}
\put(1881,398){\ellipse{150}{150}}
\put(1890,97){\shade\ellipse{150}{150}}
\put(1890,97){\ellipse{150}{150}}
\put(3325,239){\shade\ellipse{150}{150}}
\put(3325,239){\ellipse{150}{150}}
\put(3069,386){\shade\ellipse{150}{150}}
\put(3069,386){\ellipse{150}{150}}
\put(3076,83){\shade\ellipse{150}{150}}
\put(3076,83){\ellipse{150}{150}}
\put(3622,238){\shade\ellipse{150}{150}}
\put(3622,238){\ellipse{150}{150}}
\put(2458,168){\makebox(0,0)[lb]{\smash{{\SetFigFontNFSS{12}{14.4}{\rmdefault}{\mddefault}{\updefault}$+\,\,\frac{1}{2}$}}}}
\put(1485,174){\makebox(0,0)[lb]{\smash{{\SetFigFontNFSS{12}{14.4}{\rmdefault}{\mddefault}{\updefault}$\frac{1}{8}$}}}}
\put(15,174){\makebox(0,0)[lb]{\smash{{\SetFigFontNFSS{12}{14.4}{\rmdefault}{\mddefault}{\updefault}$\vertex=4\,,\lp=1$}}}}
\end{picture}
}

%% file: figures/n6k0.eepic
\setlength{\unitlength}{0.00083333in}
\begingroup\makeatletter\ifx\SetFigFontNFSS\undefined%
\gdef\SetFigFontNFSS#1#2#3#4#5{%
  \reset@font\fontsize{#1}{#2pt}%
  \fontfamily{#3}\fontseries{#4}\fontshape{#5}%
  \selectfont}%
\fi\endgroup%
{\renewcommand{\dashlinestretch}{30}
\begin{picture}(6528,2367)(0,-10)
\path(1810,1756)(3310,1756)
\path(2110,383)(2995,383)
\texture{44555555 55aaaaaa aa555555 55aaaaaa aa555555 55aaaaaa aa555555 55aaaaaa 
	aa555555 55aaaaaa aa555555 55aaaaaa aa555555 55aaaaaa aa555555 55aaaaaa 
	aa555555 55aaaaaa aa555555 55aaaaaa aa555555 55aaaaaa aa555555 55aaaaaa 
	aa555555 55aaaaaa aa555555 55aaaaaa aa555555 55aaaaaa aa555555 55aaaaaa }
\path(4487,1763)(4187,1763)
\path(4487,1763)(4187,1763)
\path(4037,2023)(4197,1753)
\path(4047,1513)(4207,1783)
\path(4757,1763)(4457,1763)
\path(4757,1763)(4457,1763)
\path(5087,1758)(4787,1758)
\path(5087,1758)(4787,1758)
\path(3752,658)(3912,388)
\path(3762,148)(3922,418)
\path(3917,391)(4172,391)(4217,406)
\path(4340,140)(4180,410)
\path(4330,650)(4170,380)
\path(5359,381)(5659,381)
\path(5359,381)(5659,381)
\path(5454,671)(5354,341)
\path(5454,91)(5364,391)
\path(5354,391)(5114,191)
\path(5354,371)(5114,551)
\path(6145,1748)(5845,1223)
\path(6005,1498)(6165,1768)
\path(6445,1748)(6145,1748)
\path(6445,1748)(6145,1748)
\path(5995,2008)(6155,1738)
\path(5845,2273)(6145,1748)
\path(2395,83)(2395,383)
\path(2395,83)(2395,383)
\path(2395,383)(2395,683)
\path(2395,383)(2395,683)
\path(2395,83)(2395,383)
\path(2395,83)(2395,383)
\put(4187,1763){\shade\ellipse{150}{150}}
\put(4187,1763){\ellipse{150}{150}}
\put(4487,1763){\shade\ellipse{150}{150}}
\put(4487,1763){\ellipse{150}{150}}
\put(4787,1763){\shade\ellipse{150}{150}}
\put(4787,1763){\ellipse{150}{150}}
\put(4042,2023){\shade\ellipse{150}{150}}
\put(4042,2023){\ellipse{150}{150}}
\put(4042,1503){\shade\ellipse{150}{150}}
\put(4042,1503){\ellipse{150}{150}}
\put(5092,1758){\shade\ellipse{150}{150}}
\put(5092,1758){\ellipse{150}{150}}
\put(3902,398){\shade\ellipse{150}{150}}
\put(3902,398){\ellipse{150}{150}}
\put(3757,658){\shade\ellipse{150}{150}}
\put(3757,658){\ellipse{150}{150}}
\put(3757,138){\shade\ellipse{150}{150}}
\put(3757,138){\ellipse{150}{150}}
\put(4190,400){\shade\ellipse{150}{150}}
\put(4190,400){\ellipse{150}{150}}
\put(4335,140){\shade\ellipse{150}{150}}
\put(4335,140){\ellipse{150}{150}}
\put(4335,660){\shade\ellipse{150}{150}}
\put(4335,660){\ellipse{150}{150}}
\put(5359,381){\shade\ellipse{150}{150}}
\put(5359,381){\ellipse{150}{150}}
\put(5659,381){\shade\ellipse{150}{150}}
\put(5659,381){\ellipse{150}{150}}
\put(5454,86){\shade\ellipse{150}{150}}
\put(5454,86){\ellipse{150}{150}}
\put(5124,186){\shade\ellipse{150}{150}}
\put(5124,186){\ellipse{150}{150}}
\put(5114,556){\shade\ellipse{150}{150}}
\put(5114,556){\ellipse{150}{150}}
\put(5454,666){\shade\ellipse{150}{150}}
\put(5454,666){\ellipse{150}{150}}
\put(6000,1488){\shade\ellipse{150}{150}}
\put(6000,1488){\ellipse{150}{150}}
\put(5854,1231){\shade\ellipse{150}{150}}
\put(5854,1231){\ellipse{150}{150}}
\put(6145,1748){\shade\ellipse{150}{150}}
\put(6145,1748){\ellipse{150}{150}}
\put(6445,1748){\shade\ellipse{150}{150}}
\put(6445,1748){\ellipse{150}{150}}
\put(6000,2008){\shade\ellipse{150}{150}}
\put(6000,2008){\ellipse{150}{150}}
\put(5850,2269){\shade\ellipse{150}{150}}
\put(5850,2269){\ellipse{150}{150}}
\put(2095,383){\shade\ellipse{150}{150}}
\put(2095,383){\ellipse{150}{150}}
\put(2395,83){\shade\ellipse{150}{150}}
\put(2395,83){\ellipse{150}{150}}
\put(2395,683){\shade\ellipse{150}{150}}
\put(2395,683){\ellipse{150}{150}}
\put(2695,383){\shade\ellipse{150}{150}}
\put(2695,383){\ellipse{150}{150}}
\put(2980,376){\shade\ellipse{150}{150}}
\put(2980,376){\ellipse{150}{150}}
\put(1802,1755){\shade\ellipse{150}{150}}
\put(1802,1755){\ellipse{150}{150}}
\put(2102,1755){\shade\ellipse{150}{150}}
\put(2102,1755){\ellipse{150}{150}}
\put(2402,1755){\shade\ellipse{150}{150}}
\put(2402,1755){\ellipse{150}{150}}
\put(3002,1755){\shade\ellipse{150}{150}}
\put(3002,1755){\ellipse{150}{150}}
\put(3302,1755){\shade\ellipse{150}{150}}
\put(3302,1755){\ellipse{150}{150}}
\put(2702,1755){\shade\ellipse{150}{150}}
\put(2702,1755){\ellipse{150}{150}}
\put(2387,383){\shade\ellipse{150}{150}}
\put(2387,383){\ellipse{150}{150}}
\put(15,1691){\makebox(0,0)[lb]{\smash{{\SetFigFontNFSS{12}{14.4}{\rmdefault}{\mddefault}{\updefault}$\vertex=6\,,\lp=0$}}}}
\put(1435,1701){\makebox(0,0)[lb]{\smash{{\SetFigFontNFSS{12}{14.4}{\rmdefault}{\mddefault}{\updefault}$\frac{1}{2}$}}}}
\put(3472,1678){\makebox(0,0)[lb]{\smash{{\SetFigFontNFSS{12}{14.4}{\rmdefault}{\mddefault}{\updefault}$+\,\,\frac{1}{2}$}}}}
\put(5282,1691){\makebox(0,0)[lb]{\smash{{\SetFigFontNFSS{12}{14.4}{\rmdefault}{\mddefault}{\updefault}$+\,\,\frac{1}{2}$}}}}
\put(1512,341){\makebox(0,0)[lb]{\smash{{\SetFigFontNFSS{12}{14.4}{\rmdefault}{\mddefault}{\updefault}$+\,\,\frac{1}{3!}$}}}}
\put(3190,341){\makebox(0,0)[lb]{\smash{{\SetFigFontNFSS{12}{14.4}{\rmdefault}{\mddefault}{\updefault}$+\,\,\frac{1}{2^3}$}}}}
\put(4570,351){\makebox(0,0)[lb]{\smash{{\SetFigFontNFSS{12}{14.4}{\rmdefault}{\mddefault}{\updefault}$+\,\,\frac{1}{5!}$}}}}
\end{picture}
}

%% file: figures/simplen5k1.eepic
\setlength{\unitlength}{0.00083333in}
\begingroup\makeatletter\ifx\SetFigFontNFSS\undefined%
\gdef\SetFigFontNFSS#1#2#3#4#5{%
  \reset@font\fontsize{#1}{#2pt}%
  \fontfamily{#3}\fontseries{#4}\fontshape{#5}%
  \selectfont}%
\fi\endgroup%
{\renewcommand{\dashlinestretch}{30}
\begin{picture}(6569,701)(0,-10)
\path(2169,500)(1771,500)
\path(2138,205)(1740,206)
\path(3290,605)(3572,324)
\path(3070,393)(3351,112)
\path(2102,500)(2416,360)
\path(2416,353)(2102,193)
\texture{44555555 55aaaaaa aa555555 55aaaaaa aa555555 55aaaaaa aa555555 55aaaaaa 
	aa555555 55aaaaaa aa555555 55aaaaaa aa555555 55aaaaaa aa555555 55aaaaaa 
	aa555555 55aaaaaa aa555555 55aaaaaa aa555555 55aaaaaa aa555555 55aaaaaa 
	aa555555 55aaaaaa aa555555 55aaaaaa aa555555 55aaaaaa aa555555 55aaaaaa }
\path(1788,493)(1788,194)
\path(1788,493)(1788,194)
\path(3326,162)(3538,374)
\path(3326,162)(3538,374)
\path(3113,375)(3326,587)
\path(3113,375)(3326,587)
\path(3535,370)(3835,370)
\path(3535,370)(3835,370)
\path(4828,349)(4580,492)(4580,206)(4828,349)
\path(4901,353)(5201,353)
\path(4901,353)(5201,353)
\path(5174,346)(5474,346)
\path(5174,346)(5474,346)
\path(6238,339)(6486,196)(6486,482)(6238,339)
\path(6081,603)(6241,333)
\path(6091,93)(6251,363)
\put(1800,206){\shade\ellipse{150}{150}}
\put(1800,206){\ellipse{150}{150}}
\put(1800,498){\shade\ellipse{150}{150}}
\put(1800,498){\ellipse{150}{150}}
\put(2095,500){\shade\ellipse{150}{150}}
\put(2095,500){\ellipse{150}{150}}
\put(2402,359){\shade\ellipse{150}{150}}
\put(2402,359){\ellipse{150}{150}}
\put(2094,202){\shade\ellipse{150}{150}}
\put(2094,202){\ellipse{150}{150}}
\put(3320,574){\shade\ellipse{150}{150}}
\put(3320,574){\ellipse{150}{150}}
\put(3528,366){\shade\ellipse{150}{150}}
\put(3528,366){\ellipse{150}{150}}
\put(3100,360){\shade\ellipse{150}{150}}
\put(3100,360){\ellipse{150}{150}}
\put(3319,153){\shade\ellipse{150}{150}}
\put(3319,153){\ellipse{150}{150}}
\put(3815,369){\shade\ellipse{150}{150}}
\put(3815,369){\ellipse{150}{150}}
\put(4836,352){\shade\ellipse{150}{150}}
\put(4836,352){\ellipse{150}{150}}
\put(4580,499){\shade\ellipse{150}{150}}
\put(4580,499){\ellipse{150}{150}}
\put(4587,196){\shade\ellipse{150}{150}}
\put(4587,196){\ellipse{150}{150}}
\put(5133,351){\shade\ellipse{150}{150}}
\put(5133,351){\ellipse{150}{150}}
\put(5426,351){\shade\ellipse{150}{150}}
\put(5426,351){\ellipse{150}{150}}
\put(6230,336){\shade\ellipse{150}{150}}
\put(6230,336){\ellipse{150}{150}}
\put(6486,189){\shade\ellipse{150}{150}}
\put(6486,189){\ellipse{150}{150}}
\put(6479,492){\shade\ellipse{150}{150}}
\put(6479,492){\ellipse{150}{150}}
\put(6086,603){\shade\ellipse{150}{150}}
\put(6086,603){\ellipse{150}{150}}
\put(6086,83){\shade\ellipse{150}{150}}
\put(6086,83){\ellipse{150}{150}}
\put(15,298){\makebox(0,0)[lb]{\smash{{\SetFigFontNFSS{12}{14.4}{\rmdefault}{\mddefault}{\updefault}$\vertex=5\,,\lp=1$}}}}
\put(1436,298){\makebox(0,0)[lb]{\smash{{\SetFigFontNFSS{12}{14.4}{\rmdefault}{\mddefault}{\updefault}$\frac{1}{10}$}}}}
\put(2638,301){\makebox(0,0)[lb]{\smash{{\SetFigFontNFSS{12}{14.4}{\rmdefault}{\mddefault}{\updefault}$+\,\,\frac{1}{2}$}}}}
\put(4018,287){\makebox(0,0)[lb]{\smash{{\SetFigFontNFSS{12}{14.4}{\rmdefault}{\mddefault}{\updefault}$+\,\,\frac{1}{2}$}}}}
\put(5614,291){\makebox(0,0)[lb]{\smash{{\SetFigFontNFSS{12}{14.4}{\rmdefault}{\mddefault}{\updefault}$+\,\,\frac{1}{4}$}}}}
\end{picture}
}

%% file: figures/simplen4k2.eepic
\setlength{\unitlength}{0.00083333in}
\begingroup\makeatletter\ifx\SetFigFontNFSS\undefined%
\gdef\SetFigFontNFSS#1#2#3#4#5{%
  \reset@font\fontsize{#1}{#2pt}%
  \fontfamily{#3}\fontseries{#4}\fontshape{#5}%
  \selectfont}%
\fi\endgroup%
{\renewcommand{\dashlinestretch}{30}
\begin{picture}(2181,482)(0,-10)
\path(2099,423)(2099,25)
\path(1793,429)(1793,31)
\path(1795,393)(2085,103)(2095,113)
\texture{44555555 55aaaaaa aa555555 55aaaaaa aa555555 55aaaaaa aa555555 55aaaaaa 
	aa555555 55aaaaaa aa555555 55aaaaaa aa555555 55aaaaaa aa555555 55aaaaaa 
	aa555555 55aaaaaa aa555555 55aaaaaa aa555555 55aaaaaa aa555555 55aaaaaa 
	aa555555 55aaaaaa aa555555 55aaaaaa aa555555 55aaaaaa aa555555 55aaaaaa }
\path(1811,84)(2111,84)
\path(1811,84)(2111,84)
\path(1811,385)(2111,385)
\path(1811,385)(2111,385)
\put(2098,380){\shade\ellipse{150}{150}}
\put(2098,380){\ellipse{150}{150}}
\put(2098,85){\shade\ellipse{150}{150}}
\put(2098,85){\ellipse{150}{150}}
\put(1791,384){\shade\ellipse{150}{150}}
\put(1791,384){\ellipse{150}{150}}
\put(1800,83){\shade\ellipse{150}{150}}
\put(1800,83){\ellipse{150}{150}}
\put(1415,170){\makebox(0,0)[lb]{\smash{{\SetFigFontNFSS{12}{14.4}{\rmdefault}{\mddefault}{\updefault}$\frac{1}{4}$}}}}
\put(15,170){\makebox(0,0)[lb]{\smash{{\SetFigFontNFSS{12}{14.4}{\rmdefault}{\mddefault}{\updefault}$\vertex=4\,,\lp=2$}}}}
\end{picture}
}

%% file: appnoloops.tex
\section{}\label{app:noloops}
This appendix shows  $\Ds^{\vertex,\lp,0}$ or $\Dsnew^{\vertex,\lp,0}$ up to order $\vertex+\lp\le5$ and computed via formulas  (\ref{eq:recnoloops}) or  (\ref{eq:rectheta}), respectively. All graphs in the same equivalence class are identified as the same.   The coefficients in front  of graphs are the
inverses of the orders of their
groups of automorphisms.
\vspace{1cm}\\
\input{figures/n1k0.eepic}
\vspace{0.95cm}\\
\input{figures/n2k0.eepic}
\vspace{0.95cm}\\
\input{figures/n3k0.eepic}
\vspace{0.9cm}\\
\input{figures/noloopsn2k1.eepic}
\vspace{0.5cm}\\
\input{figures/n4k0.eepic}
\vspace{0.9cm}\\
\input{figures/noloopsn3k1.eepic}
\vspace{1cm}\\
\input{figures/noloopsn2k2.eepic}
\vspace{0.5cm}\\
\input{figures/n5k0.eepic}
\vspace{0.9cm}\\
\input{figures/noloopsn4k1.eepic}
\vspace{1cm}\\
\input{figures/noloopsn3k2.eepic}
\vspace{1cm}\\
\input{figures/noloopsn2k3.eepic}

%% file: figures/noloopsn3k1.eepic
\setlength{\unitlength}{0.00083333in}
\begingroup\makeatletter\ifx\SetFigFontNFSS\undefined%
\gdef\SetFigFontNFSS#1#2#3#4#5{%
  \reset@font\fontsize{#1}{#2pt}%
  \fontfamily{#3}\fontseries{#4}\fontshape{#5}%
  \selectfont}%
\fi\endgroup%
{\renewcommand{\dashlinestretch}{30}
\begin{picture}(3796,435)(0,-10)
\put(3280,160){\ellipse{300}{150}}
\path(3389,157)(3749,157)(3755,157)
\path(2143,331)(2000,83)(2286,83)(2143,331)
\texture{44555555 55aaaaaa aa555555 55aaaaaa aa555555 55aaaaaa aa555555 55aaaaaa 
	aa555555 55aaaaaa aa555555 55aaaaaa aa555555 55aaaaaa aa555555 55aaaaaa 
	aa555555 55aaaaaa aa555555 55aaaaaa aa555555 55aaaaaa aa555555 55aaaaaa 
	aa555555 55aaaaaa aa555555 55aaaaaa aa555555 55aaaaaa aa555555 55aaaaaa }
\put(2140,339){\shade\ellipse{150}{150}}
\put(2140,339){\ellipse{150}{150}}
\put(1993,83){\shade\ellipse{150}{150}}
\put(1993,83){\ellipse{150}{150}}
\put(2296,90){\shade\ellipse{150}{150}}
\put(2296,90){\ellipse{150}{150}}
\put(3131,160){\shade\ellipse{150}{150}}
\put(3131,160){\ellipse{150}{150}}
\put(3425,160){\shade\ellipse{150}{150}}
\put(3425,160){\ellipse{150}{150}}
\put(3713,166){\shade\ellipse{150}{150}}
\put(3713,166){\ellipse{150}{150}}
\put(1515,80){\makebox(0,0)[lb]{\smash{{\SetFigFontNFSS{12}{14.4}{\rmdefault}{\mddefault}{\updefault}$\frac{1}{3!}$}}}}
\put(15,95){\makebox(0,0)[lb]{\smash{{\SetFigFontNFSS{12}{14.4}{\rmdefault}{\mddefault}{\updefault}$\lp=1,\vertex=3$}}}}
\put(2555,90){\makebox(0,0)[lb]{\smash{{\SetFigFontNFSS{12}{14.4}{\rmdefault}{\mddefault}{\updefault}$+\,\,\frac{1}{2}$}}}}
\end{picture}
}

%% file: figures/noloopsn2k2.eepic
\setlength{\unitlength}{0.00083333in}
\begingroup\makeatletter\ifx\SetFigFontNFSS\undefined%
\gdef\SetFigFontNFSS#1#2#3#4#5{%
  \reset@font\fontsize{#1}{#2pt}%
  \fontfamily{#3}\fontseries{#4}\fontshape{#5}%
  \selectfont}%
\fi\endgroup%
{\renewcommand{\dashlinestretch}{30}
\begin{picture}(2285,276)(0,-10)
\put(2056,144){\ellipse{300}{150}}
\path(1878,141)(2202,141)(2190,141)
\texture{44555555 55aaaaaa aa555555 55aaaaaa aa555555 55aaaaaa aa555555 55aaaaaa 
	aa555555 55aaaaaa aa555555 55aaaaaa aa555555 55aaaaaa aa555555 55aaaaaa 
	aa555555 55aaaaaa aa555555 55aaaaaa aa555555 55aaaaaa aa555555 55aaaaaa 
	aa555555 55aaaaaa aa555555 55aaaaaa aa555555 55aaaaaa aa555555 55aaaaaa }
\put(2202,144){\shade\ellipse{150}{150}}
\put(2202,144){\ellipse{150}{150}}
\put(1908,144){\shade\ellipse{150}{150}}
\put(1908,144){\ellipse{150}{150}}
\put(15,78){\makebox(0,0)[lb]{\smash{{\SetFigFontNFSS{12}{14.4}{\rmdefault}{\mddefault}{\updefault}$\vertex=2\,,\lp=2$}}}}
\put(1440,78){\makebox(0,0)[lb]{\smash{{\SetFigFontNFSS{12}{14.4}{\rmdefault}{\mddefault}{\updefault}$\frac{1}{2\cdot3!}$}}}}
\end{picture}
}

%% file: figures/noloopsn4k1.eepic
\setlength{\unitlength}{0.00083333in}
\begingroup\makeatletter\ifx\SetFigFontNFSS\undefined%
\gdef\SetFigFontNFSS#1#2#3#4#5{%
  \reset@font\fontsize{#1}{#2pt}%
  \fontfamily{#3}\fontseries{#4}\fontshape{#5}%
  \selectfont}%
\fi\endgroup%
{\renewcommand{\dashlinestretch}{30}
\begin{picture}(6460,1212)(0,-10)
\path(2115,1120)(2115,722)
\path(1809,1126)(1809,728)
\path(3595,909)(3235,909)(3229,909)
\put(5641,875){\ellipse{300}{150}}
\path(5750,872)(6110,872)(6116,872)
\path(6445,863)(6085,863)(6079,863)
\put(2430,145){\ellipse{300}{150}}
\path(2539,142)(2899,142)(2905,142)
\path(2334,146)(1974,146)(1968,146)
\put(4629,864){\ellipse{300}{150}}
\texture{44555555 55aaaaaa aa555555 55aaaaaa aa555555 55aaaaaa aa555555 55aaaaaa 
	aa555555 55aaaaaa aa555555 55aaaaaa aa555555 55aaaaaa aa555555 55aaaaaa 
	aa555555 55aaaaaa aa555555 55aaaaaa aa555555 55aaaaaa aa555555 55aaaaaa 
	aa555555 55aaaaaa aa555555 55aaaaaa aa555555 55aaaaaa aa555555 55aaaaaa }
\path(1827,781)(2127,781)
\path(1827,781)(2127,781)
\path(1827,1082)(2127,1082)
\path(1827,1082)(2127,1082)
\path(3219,907)(2971,1050)(2971,764)(3219,907)
\path(4351,1114)(4511,844)
\path(4361,604)(4521,874)
\put(2114,1077){\shade\ellipse{150}{150}}
\put(2114,1077){\ellipse{150}{150}}
\put(2114,782){\shade\ellipse{150}{150}}
\put(2114,782){\ellipse{150}{150}}
\put(1807,1081){\shade\ellipse{150}{150}}
\put(1807,1081){\ellipse{150}{150}}
\put(1816,780){\shade\ellipse{150}{150}}
\put(1816,780){\ellipse{150}{150}}
\put(3227,910){\shade\ellipse{150}{150}}
\put(3227,910){\ellipse{150}{150}}
\put(2971,1057){\shade\ellipse{150}{150}}
\put(2971,1057){\ellipse{150}{150}}
\put(2978,754){\shade\ellipse{150}{150}}
\put(2978,754){\ellipse{150}{150}}
\put(3522,907){\shade\ellipse{150}{150}}
\put(3522,907){\ellipse{150}{150}}
\put(5492,875){\shade\ellipse{150}{150}}
\put(5492,875){\ellipse{150}{150}}
\put(5786,875){\shade\ellipse{150}{150}}
\put(5786,875){\ellipse{150}{150}}
\put(6074,881){\shade\ellipse{150}{150}}
\put(6074,881){\ellipse{150}{150}}
\put(6377,885){\shade\ellipse{150}{150}}
\put(6377,885){\ellipse{150}{150}}
\put(2281,145){\shade\ellipse{150}{150}}
\put(2281,145){\ellipse{150}{150}}
\put(2575,145){\shade\ellipse{150}{150}}
\put(2575,145){\ellipse{150}{150}}
\put(2863,151){\shade\ellipse{150}{150}}
\put(2863,151){\ellipse{150}{150}}
\put(2003,142){\shade\ellipse{150}{150}}
\put(2003,142){\ellipse{150}{150}}
\put(4480,864){\shade\ellipse{150}{150}}
\put(4480,864){\ellipse{150}{150}}
\put(4774,864){\shade\ellipse{150}{150}}
\put(4774,864){\ellipse{150}{150}}
\put(4356,1114){\shade\ellipse{150}{150}}
\put(4356,1114){\ellipse{150}{150}}
\put(4356,594){\shade\ellipse{150}{150}}
\put(4356,594){\ellipse{150}{150}}
\put(15,815){\makebox(0,0)[lb]{\smash{{\SetFigFontNFSS{12}{14.4}{\rmdefault}{\mddefault}{\updefault}$\vertex=4\,,\lp=1$}}}}
\put(1431,817){\makebox(0,0)[lb]{\smash{{\SetFigFontNFSS{12}{14.4}{\rmdefault}{\mddefault}{\updefault}$\frac{1}{8}$}}}}
\put(2352,813){\makebox(0,0)[lb]{\smash{{\SetFigFontNFSS{12}{14.4}{\rmdefault}{\mddefault}{\updefault}$+\,\,\frac{1}{2}$}}}}
\put(3751,817){\makebox(0,0)[lb]{\smash{{\SetFigFontNFSS{12}{14.4}{\rmdefault}{\mddefault}{\updefault}$+\,\,\frac{1}{4}$}}}}
\put(4970,837){\makebox(0,0)[lb]{\smash{{\SetFigFontNFSS{12}{14.4}{\rmdefault}{\mddefault}{\updefault}$+\,\,\frac{1}{2}$}}}}
\put(1476,78){\makebox(0,0)[lb]{\smash{{\SetFigFontNFSS{12}{14.4}{\rmdefault}{\mddefault}{\updefault}$+\,\,\frac{1}{4}$}}}}
\end{picture}
}

%% file: figures/noloopsn3k2.eepic
\setlength{\unitlength}{0.00083333in}
\begingroup\makeatletter\ifx\SetFigFontNFSS\undefined%
\gdef\SetFigFontNFSS#1#2#3#4#5{%
  \reset@font\fontsize{#1}{#2pt}%
  \fontfamily{#3}\fontseries{#4}\fontshape{#5}%
  \selectfont}%
\fi\endgroup%
{\renewcommand{\dashlinestretch}{30}
\begin{picture}(5403,549)(0,-10)
\path(1946,84)(2093,457)
\path(2253,90)(2099,457)
\drawline(2099,457)(2099,457)
\put(2106,84){\ellipse{300}{150}}
\put(3388,233){\ellipse{300}{150}}
\put(3690,227){\ellipse{300}{150}}
\put(4864,221){\ellipse{300}{150}}
\path(4686,218)(5010,218)(4998,218)
\path(5002,220)(5326,220)(5314,220)
\texture{44555555 55aaaaaa aa555555 55aaaaaa aa555555 55aaaaaa aa555555 55aaaaaa 
	aa555555 55aaaaaa aa555555 55aaaaaa aa555555 55aaaaaa aa555555 55aaaaaa 
	aa555555 55aaaaaa aa555555 55aaaaaa aa555555 55aaaaaa aa555555 55aaaaaa 
	aa555555 55aaaaaa aa555555 55aaaaaa aa555555 55aaaaaa aa555555 55aaaaaa }
\put(1950,81){\shade\ellipse{150}{150}}
\put(1950,81){\ellipse{150}{150}}
\put(2250,87){\shade\ellipse{150}{150}}
\put(2250,87){\ellipse{150}{150}}
\put(2097,453){\shade\ellipse{150}{150}}
\put(2097,453){\ellipse{150}{150}}
\put(3239,233){\shade\ellipse{150}{150}}
\put(3239,233){\ellipse{150}{150}}
\put(3533,233){\shade\ellipse{150}{150}}
\put(3533,233){\ellipse{150}{150}}
\put(3838,234){\shade\ellipse{150}{150}}
\put(3838,234){\ellipse{150}{150}}
\put(5010,221){\shade\ellipse{150}{150}}
\put(5010,221){\ellipse{150}{150}}
\put(4716,221){\shade\ellipse{150}{150}}
\put(4716,221){\ellipse{150}{150}}
\put(5320,211){\shade\ellipse{150}{150}}
\put(5320,211){\ellipse{150}{150}}
\put(15,137){\makebox(0,0)[lb]{\smash{{\SetFigFontNFSS{12}{14.4}{\rmdefault}{\mddefault}{\updefault}$\vertex=3\,,\lp=2$}}}}
\put(1446,126){\makebox(0,0)[lb]{\smash{{\SetFigFontNFSS{12}{14.4}{\rmdefault}{\mddefault}{\updefault}$\frac{1}{2^2}$}}}}
\put(2641,151){\makebox(0,0)[lb]{\smash{{\SetFigFontNFSS{12}{14.4}{\rmdefault}{\mddefault}{\updefault}$+\,\,\frac{1}{2^3}$}}}}
\put(4150,152){\makebox(0,0)[lb]{\smash{{\SetFigFontNFSS{12}{14.4}{\rmdefault}{\mddefault}{\updefault}$+\,\,\frac{1}{3!}$}}}}
\end{picture}
}

%% file: figures/noloopsn2k3.eepic
\setlength{\unitlength}{0.00083333in}
\begingroup\makeatletter\ifx\SetFigFontNFSS\undefined%
\gdef\SetFigFontNFSS#1#2#3#4#5{%
  \reset@font\fontsize{#1}{#2pt}%
  \fontfamily{#3}\fontseries{#4}\fontshape{#5}%
  \selectfont}%
\fi\endgroup%
{\renewcommand{\dashlinestretch}{30}
\begin{picture}(2343,341)(0,-10)
\put(2117,163){\ellipse{310}{310}}
\put(2112,161){\ellipse{300}{150}}
\texture{44555555 55aaaaaa aa555555 55aaaaaa aa555555 55aaaaaa aa555555 55aaaaaa 
	aa555555 55aaaaaa aa555555 55aaaaaa aa555555 55aaaaaa aa555555 55aaaaaa 
	aa555555 55aaaaaa aa555555 55aaaaaa aa555555 55aaaaaa aa555555 55aaaaaa 
	aa555555 55aaaaaa aa555555 55aaaaaa aa555555 55aaaaaa aa555555 55aaaaaa }
\put(1968,164){\shade\ellipse{150}{150}}
\put(1968,164){\ellipse{150}{150}}
\put(2260,167){\shade\ellipse{150}{150}}
\put(2260,167){\ellipse{150}{150}}
\put(15,92){\makebox(0,0)[lb]{\smash{{\SetFigFontNFSS{12}{14.4}{\rmdefault}{\mddefault}{\updefault}$\vertex=2\,,\lp=3$}}}}
\put(1451,95){\makebox(0,0)[lb]{\smash{{\SetFigFontNFSS{12}{14.4}{\rmdefault}{\mddefault}{\updefault}$\frac{1}{2\cdot4!}$}}}}
\end{picture}
}